\documentclass[12pt,twoside,reqno]{amsart}
\usepackage{amsmath}
\usepackage{amsfonts}
\usepackage{amssymb}
\usepackage{color}
\usepackage{mathrsfs}
\usepackage{multicol}
\usepackage{graphicx}
\usepackage{cite}
\usepackage[dvipdfmx,colorlinks=true, pdfstartview=FitV, linkcolor=blue,citecolor=blue, urlcolor=blue]{hyperref}
\usepackage{geometry}
\usepackage{marginnote}
\usepackage{todonotes}
\usepackage{tikz}
\usepackage{pgfplots}
\pgfplotsset{compat=1.10}
\usepgfplotslibrary{fillbetween}
\usetikzlibrary{patterns}
\newenvironment{pr}[1]
{{\noindent \bf Proof of {#1}.\  }}{\hfill \qed}
\DeclareMathOperator*{\supp}{supp}
\newtheorem{theorem}{Theorem}[section]

\newtheorem{lemma}[theorem]{Lemma}
\newtheorem{proposition}[theorem]{Proposition}

\newtheorem{remark}[theorem]{Remark}
\numberwithin{equation}{section}

\newcommand{\eqntag}{\addtocounter{equation}{1}\tag{\theequation}} 
\allowdisplaybreaks[1]
\begin{document}
\title{Boundedness of radial solutions to the parabolic-parabolic Keller--Segel system on the plane with subcritical mass} 
\author{Tatsuya Hosono}
\address{Osaka Central Advanced Mathematical Institute, Osaka Metropolitan University,
	Osaka 558-8585, Japan}
\email{tatsuya.hosono@omu.ac.jp}
	
	
\makeatletter
\def\subjclassname{2020 Mathematics Subject Classification}
\makeatother

\date{\today}

\begin{abstract}
We establish uniform-in-time boundedness of radially symmetric global solutions to the two-dimensional parabolic-parabolic Keller–Segel system on the whole space with mass below $8\pi$. While global existence is known throughout this subcritical mass range, uniform-in-time boundedness has remained open. On bounded domains, the free energy method of Nagai, Senba and Yoshida (1997) yields both global existence and boundedness under the optimal mass condition. On the whole space, however, the corresponding estimates are not uniform in time due to the behavior of solutions at spatial infinity. Indeed, the free energy of self-similar solutions diverges to $-\infty$ as $t\to\infty$, preventing a direct application of this approach to establish boundedness. The key ingredients are new interior estimates controlling the localized free energy in terms of its dissipation, together with uniform-in-time exterior estimates exploiting radial symmetry.
\end{abstract}
%
\maketitle
%
%
\pagestyle{myheadings}
\markboth
{\sc{Tatsuya Hosono}}
{\sc{Boundedness for the Keller--Segel system on the plane}}

\section{Introduction}
Chemotaxis is the movement of cells in response to spatial variations in the concentration of a chemical substance. In the case of an attractive signal, cells tend to migrate toward regions of higher concentrations of the chemoattractant in addition to undergoing random diffusion. When the signal is secreted by the cells themselves, this mechanism can lead to the spontaneous formation of cell aggregates. The Keller–Segel system \cite{Ke-Se70,Pa} describes this process through coupled equations for the cell density and the concentration of the chemical signal.
In this paper, we study the parabolic-parabolic Keller–Segel system on the plane
\begin{align}
\left\{
\begin{aligned}
&\partial_tu -\Delta u + \nabla \cdot (u\nabla v) =0, & t>0,~~&x\in\mathbb{R}^2,
\\
&\tau \partial _t v-\Delta v +\lambda v= u,  & t>0,~~&x\in\mathbb{R}^2,
\\
&(u,v)(0,x)=(u_0,v_0)(x),&&x\in \mathbb{R}^2,
\end{aligned}
\right.
\label{eqn;KS}
\end{align}
where $\tau,\lambda>0$ and $u_0,v_0\ge0$. The unknown functions $u=u(t,x)$ and $v=v(t,x)$ denote the cell density and the concentration of the chemoattractant, respectively. 
The first equation describes cell diffusion and migration toward higher concentrations of the chemoattractant, while the second accounts for its diffusion, production by the cells, and degradation.
Such aggregation phenomena are commonly interpreted mathematically as blowup of the cell density, in the sense that
$\lim_{t\to T}\|u(t)\|_{\infty}=\infty$
for some $T\in(0,\infty]$.
Determining conditions under which solutions remain bounded or undergo blowup is therefore a central issue in the mathematical analysis of this system.
In two space dimensions, under mass conservation $\|u(t)\|_1=\|u_0\|_1$,
 the total cell mass plays a fundamental role in the behavior of solutions. Although global existence is known throughout the subcritical mass $\|u_0\|_1<8\pi$, it does not by itself exclude unbounded growth as $t\to\infty$. The purpose of this paper is to establish uniform-in-time boundedness throughout this mass range for radially symmetric solutions, without imposing any moment condition on the initial cell density.

For the system \eqref{eqn;KS}, global existence in the subcritical regime $\|u_0\|_1<8\pi$ was established by Calvez and Corrias \cite{Ca-Co08} under additional assumptions of finite initial entropy and logarithmic moment. Nagai and Ogawa \cite{Na-Og11} subsequently obtained global existence without these assumptions under the stronger mass restriction $\|u_0\|_1<4\pi$. Mizoguchi \cite{Mi13} extended the latter result to the full subcritical range $\|u_0\|_1<8\pi$.
In the supercritical regime $\|u_0\|_1>8\pi$, finite-time blowup solutions on the plane were constructed by Mizoguchi \cite{Mi20JMPA,Mi20SIAM} in the radially symmetric setting with finite second moment. At the critical mass $\|u_0\|_1=8\pi$, global existence was established in \cite{Mi13} under radial symmetry or a logarithmic moment condition. More recently, Hosono \cite{Ho26} removed both restrictions and proved global existence for nonnegative initial data satisfying
$(u_0,v_0)\in L^1(\mathbb{R}^2)\times\bigl(L^1(\mathbb{R}^2)\cap\dot H^1(\mathbb{R}^2)\bigr)$
and $\|u_0\|_1=8\pi$. For positive forward self-similar solutions in the case $\lambda=0$, we refer to \cite{Bi-Co-Do}.
Global existence for \eqref{eqn;KS} with sufficiently large $\tau$ is also investigated in \cite{Bi-Bo-Br,Bi-Gu-Ka}.
For the corresponding Cauchy--Neumann problem in bounded domains, we refer to
\cite{Bi98,He-Ve,Ho-Wa,Na-Se-Yo}, for instance.

As mentioned before, global existence, however, does not by itself imply uniform-in-time boundedness. Even when finite-time blowup is excluded, the available estimates do not necessarily rule out unbounded growth of the cell density as $t\to\infty$. This distinction is particularly important on the whole space. For the corresponding Cauchy--Neumann problem on bounded domains $\Omega\subset\mathbb{R}^2$, Nagai, Senba and Yoshida \cite{Na-Se-Yo} obtained both global existence and uniform-in-time boundedness by combining the free energy structure with the Trudinger--Moser inequality and its optimal constant. More precisely, the free energy
\begin{align*}
\mathcal{F}(t)
=
\int_{\Omega}u\ln u\,\mathrm{d}x
-\int_{\Omega}uv\,\mathrm{d}x
+\frac12\|\nabla v\|_{L^2(\Omega)}^2
+\frac{\lambda}{2}\|v\|_{L^2(\Omega)}^2
\end{align*}
satisfies
\begin{align*}
\frac{\mathrm{d}}{\mathrm{d}t}\mathcal{F}(t)
+
\int_{\Omega}u|\nabla(\ln u-v)|^2\,\mathrm{d}x
+
\tau\|\partial_t v\|_{L^2(\Omega)}^2
=0.
\end{align*}
Together with the Trudinger--Moser inequality, this identity yields the entropy estimates needed to establish boundedness under the mass condition $\|u_0\|_1<4\pi$ in general domains and $\|u_0\|_1<8\pi$ for radially symmetric solutions in a disk.
On the whole space, a direct application of this approach encounters an additional difficulty caused by the behavior of solutions as $|x|\to\infty$. Unlike on bounded domains, the negative part of the entropy cannot be controlled by the mass alone. A logarithmic moment condition provides the estimate
\begin{align*}
\int_{\{u<1\}}u|\ln u|\,\mathrm{d}x
\le
C\int_{\mathbb{R}^2}u\ln(1+|x|^2)\,\mathrm{d}x
+C(\|u_0\|_1),
\end{align*}
but the resulting bound may grow  in time, see \cite[Theorem 3.2]{Ca-Co08}. Such growth reflects the spreading of diffusive solutions: even for the linear heat equation, the entropy of the Gaussian solution $G_t$ tends to $-\infty$ logarithmically as $t\to\infty$, i.e.,
\begin{align*}
\int_{ \mathbb R^2}G_t \ln G_t \,\mathrm{d}x=-\ln(4\pi t) -1, \qquad G_t(x)=(4\pi t)^{-1} \exp\left(-\frac{|x|^2}{4t}\right),\qquad x\in\mathbb{R}^2.
\end{align*} 
Indeed, the free energy of all self-similar solutions also diverges to $-\infty$ as $t\to\infty$, in contrast to the corresponding Cauchy--Neumann problem on bounded domains, see \cite{Mi20JMPA}. 
Thus, on the whole space, the classical entropy argument does not directly yield uniform-in-time boundedness, in contrast to the corresponding problem on bounded domains.
Alternatively, the modified entropy
\begin{align*}
\int_{\mathbb{R}^2}(1+u)\ln(1+u)\,\mathrm{d}x,
\end{align*}
introduced by Nagai \cite{Na11}, allows one to dispense with moment assumptions. However, the associated energy estimates used to establish global existence do not directly provide uniform-in-time bounds.

Boundedness on the whole space has been obtained under more restrictive assumptions. Nagai, Syukuinn and Umesako \cite{Na-Sy-Um} established boundedness for \eqref{eqn;KS} under the smallness conditions $\|u_0\|_1<\delta$ and $\|\nabla v_0\|_{L^1\cap L^\infty}<\delta$ for some $\delta>0$.
Subsequently, Nagai, Seki and Yamada \cite{Na-Se-Ya} obtained boundedness under the mass condition $\|u_0\|_1<4\pi$ in the setting of an attraction-repulsion chemotaxis system, see also \cite{Na-Ya20}. For the system with $\lambda=0$, convergence toward self-similar solutions has also been studied under additional smallness assumptions on $\tau$ and on the initial data, see \cite{Ca-Mi,Co-Es-Ma,Ka-Ma,Nai}. Nevertheless, boundedness throughout the full subcritical mass range $\|u_0\|_1<8\pi$ has remained open.
For the parabolic-elliptic Keller--Segel system with $\tau=0$,
boundedness has been established in~\cite{Di-Na-Ra,Na11,Na-Ya20}
by exploiting comparison principles in conjunction with
radial symmetry or rearrangement arguments.
However, the application of similar comparison principles
to the parabolic-parabolic system is more challenging
due to the parabolic nature of the second equation in \eqref{eqn;KS}.

In this paper, we establish uniform-in-time boundedness throughout the full subcritical mass range for radially symmetric solutions, without imposing any moment condition on the initial cell density. 
Our main result reads as follows.

\begin{theorem}\label{thm}
	Let $\tau,\lambda>0$, and let
	$(u_0,v_0)\in L^1_+(\mathbb{R}^2)\times H^1_+(\mathbb{R}^2)$
	be radially symmetric. Suppose that
\begin{align*}
\|u_0\|_1<8\pi.
\eqntag
\label{eqn;mass-constraint}
\end{align*}
Then, the radially symmetric nonnegative global solution $(u,v)$ to \eqref{eqn;KS} is bounded uniformly in time, i.e., for any $t_0>0$,
\begin{align*}
\sup_{t \ge t_0} \left(\|u(t)\|_\infty + \|\nabla v(t)\|_\infty \right) <\infty.
\end{align*}
\end{theorem}
As for the local-in-time solutions to \eqref{eqn;KS}, see for instance \cite[Proposition 3.1]{Ho26}. 
Theorem \ref{thm} establishes uniform-in-time boundedness
for every subcritical mass in the radially symmetric setting.
No smallness condition is imposed on $\tau$ or on $v_0$,
and no moment assumption is required on $u_0$.
The estimates are stated for $t\ge t_0>0$,
in accordance with the low regularity of the initial data.
Once uniform-in-time boundedness has been established, the results of \cite{Na-Sy-Um} yield the following decay estimates and asymptotic behavior for $\lambda>0$: for every $1\le p\le\infty$,
\begin{align*}
\sup_{t\ge1}t^{1-\frac1p}\|u(t)\|_p<\infty,
\qquad
\lim_{t\to\infty}t^{1-\frac1p}
\|u(t)-\|u_0\|_1G_t\|_p=0.
\end{align*}
The corresponding leading-order asymptotic profile for $v$ is
$\lambda^{-1}\|u_0\|_1G_t$, see also \cite{Na-Ya07,Ya} with $\lambda=1$.

\vspace{3mm}
\noindent
\textbf{Sketch of the proof of Theorem \ref{thm}.}
We briefly describe the main ingredients of the proof.
It is worth emphasizing that the proof of
Theorem~\ref{thm} does not rely on any comparison principle.
First, the radial estimates stated in
Appendix \ref{sect;appendix} allow us to control
the modified entropy of $u$ and the Dirichlet energy of $v$
in exterior regions.
An iteration argument then yields, for every $R>1$
and $t_0>0$,
\begin{align*}
\sup_{t\ge t_0}
\|u(t)\|_{L^\infty(\{|x|\ge R\})}<\infty.
\end{align*}
These exterior estimates do not require the mass constraint
\eqref{eqn;mass-constraint}.
The main step is then to derive a differential inequality
of the form
\begin{align*}
\frac{\mathrm{d}}{\mathrm{d}t}\mathcal{F}_R(t)
+\alpha\mathcal{F}_R(t)
+\beta\mathcal{D}_R(t)
\le C,
\eqntag
\label{eqn;ugi}
\end{align*}
where $\mathcal{F}_R$ and $\mathcal{D}_R$ denote
the localized free energy and its associated dissipation,
respectively, and $\alpha,\beta>0$,
see Proposition \ref{prop1} below.
A difficulty in deriving \eqref{eqn;ugi} is that  the dissipation
\begin{align*}
\int_{\mathbb{R}^2}
u|\nabla(\ln u-v)|^2\,\mathrm{d}x
\end{align*}
does not directly control either the standard Fisher information
\begin{align*}
\int_{ \mathbb R^2}u |\nabla \ln u|^2\,\mathrm{d}x
\end{align*}
 associated with the heat equation or the entropy due to lack of regularity of solutions.
To overcome this difficulty,
we exploit the radial structure of the system
to estimate the local entropy and Dirichlet energy
in terms of the localized dissipation:
\begin{align*}
\int_{B_R(0)}u\ln u\,\mathrm{d}x
\le\, C\bigl(1+\sqrt{\mathcal{D}_R(t)}\bigr),
\quad
\int_{B_R(0)}|\nabla v|^2\,\mathrm{d}x
\le\, C\bigl(1+\mathcal{D}_R(t)\bigr),
\end{align*}
see Lemma \ref{lem;log} below.
These are the key estimates used to establish
\eqref{eqn;ugi} and close the interior argument.
Applying Gronwall's inequality to \eqref{eqn;ugi}
then gives a uniform-in-time upper bound for $\mathcal{F}_R$.
Combining this bound with a lower estimate obtained
from the Trudinger--Moser inequality under
\eqref{eqn;mass-constraint}, we obtain
uniform-in-time interior entropy estimates.
Together with the exterior estimates,
these yield a global bound for the modified entropy,
from which the asserted boundedness follows
by standard  iteration arguments.
We expect that the approach developed in this paper
can be extended to a broader class of chemotaxis systems
on the whole space.

\vspace{3mm}
The paper is organized as follows. Section \ref{sect;2} is devoted to uniform-in-time exterior estimates without the subcritical mass restriction. In Section \ref{sect;3}, we establish the localized differential energy inequality and derive interior estimates. Section \ref{sect:4} combines these estimates to prove Theorem \ref{thm}. The radial inequalities and the annular estimate for the linear parabolic equation used in the proof are collected in Appendix \ref{sect;appendix}.

\vspace{5mm}
\section{Uniform-in-time exterior estimates}\label{sect;2}

In this section, we establish exterior estimates for solutions
using the radial lemmas stated in Appendix \ref{sect;appendix}.
The mass constraint \eqref{eqn;mass-constraint} is not required here.
The estimates developed below will be used in the next section
to derive uniform-in-time interior estimates, ultimately yielding
uniform-in-time bounds on the whole space.
Our main goal in this section is to prove that, for every~$R>1$ and~$t_0>0$,
\begin{equation*}
\sup_{t \ge t_0}
\|u(t)\|_{L^\infty(\{|x|\ge R\})}<\infty,
\end{equation*}
see Lemma \ref{lem;infty-u} below.
Together with the other exterior estimates, this bound allows us
to derive the localized differential energy inequality in
Proposition \ref{prop1}. This inequality is a key ingredient
in obtaining the interior estimates in the next section.
To this end, let $\psi\in C^\infty(\mathbb R^2)$ be a radially
symmetric cut-off function satisfying
\begin{equation}\label{eqn;cutoff}
0\le\psi\le1,
\qquad
\psi(x)=
\begin{cases}
0 & \text{if } |x|\le \dfrac1{16},\\[2mm]
1 & \text{if } |x|\ge \dfrac18.
\end{cases}
\end{equation}
For $R>1$, set $\psi_R(x):=\psi(x/R)$.
Then
\begin{align*}
|\nabla\psi_R|\le \frac{C}{R},
\qquad
|\Delta\psi_R|\le \frac{C}{R^2},
\qquad
\operatorname{supp}\nabla\psi_R
\subset
\left\{x\in\mathbb R^2;
\frac{R}{16}\le |x|\le\frac{R}{8}\right\}.
\end{align*}
We begin by establishing uniform-in-time exterior bounds
for the modified entropy of $u$ and the Dirichlet energy of $v$.
Throughout this section, we fix an arbitrary $t_0>0$ and work
on $(t_0,\infty)$, where the global solutions have the regularity
needed for the computations below.
The radius $R>1$ is arbitrary, and the constants in the estimates
may depend on $R$ and $t_0$, but are independent of time.

\begin{lemma}\label{lem;out-log}
Let $(u,v)$ be the radially symmetric solution to \eqref{eqn;KS}.
Then, for $R>1$
\begin{align*}
\sup_{t \ge t_0}\left[  \int_{\mathbb{R}^2}(1+u(t,x))\ln(1+u(t,x)) \psi_R^2(x)\,\mathrm{d}x+\int_{\mathbb{R}^2}|\nabla v(t,x)|^2\psi_R^2(x)\,\mathrm{d}x\right] <\infty
\end{align*}
and
\begin{align*}
\sup_{t \ge t_0}\int_t^{t+1} \left[ \int_{\mathbb{R}^2}\frac{ |\nabla u(s,x)|^2}{1+u(s ,x)} \psi_R^2(x)\,\mathrm{d}x+ \|\partial_tv (s)\psi_R\|_2^2 \right]\,\mathrm{d}s<\infty.
\end{align*}
\end{lemma}

\begin{pr}{Lemma \ref{lem;out-log}}
Let us first observe the time derivative of the modified entropy localized to an exterior region by means of the cut-off $\psi_R$. 
Using the first equation in \eqref{eqn;KS}, we see
\begin{align*}
&\frac{\mathrm{d}}{\mathrm{d}t}\left(\int_{\mathbb{R}^2}(1+u(t,x))\ln(1+u(t,x)) \psi_R^2(x)\,\mathrm{d}x-\int_{\mathbb{R}^2}u(t,x)\psi_R^2(x)\,\mathrm{d}x\right)
\\
=\,& -\int_{\mathbb{R}^2}\frac{|\nabla u(t,x)|^2}{1+u(t,x)} \psi_R^2(x)\,\mathrm{d}x-\int_{\mathbb{R}^2}\ln(1+u(t,x)) \nabla u(t,x)\cdot \nabla \psi_R^2(x)\,\mathrm{d}x
\\
&+\int_{\mathbb{R}^2}\frac{u(t,x)}{1+u(t,x)} \nabla u(t,x) \cdot \nabla v(t,x)\psi_R^2(x)\,\mathrm{d}x+\int_{\mathbb{R}^2}u(t,x)\ln(1+u(t,x)) \nabla v(t,x)\cdot \nabla \psi_R^2(x)\,\mathrm{d}x
\\
=\,& -\int_{\mathbb{R}^2}\frac{|\nabla u(t,x)|^2}{1+u(t,x)} \psi_R^2(x)\,\mathrm{d}x-\int_{\mathbb{R}^2}\ln(1+u(t,x)) \nabla u(t,x)\cdot \nabla \psi_R^2(x)\,\mathrm{d}x
\\
&-\int_{\mathbb{R}^2}u(t,x)\Delta v(t,x)\psi_R^2(x)\,\mathrm{d}x-\int_{\mathbb{R}^2}u(t,x) \nabla v(t,x)\cdot \nabla \psi_R^2(x)\,\mathrm{d}x
\\
&+\int_{\mathbb{R}^2}\ln(1+u(t,x)) \Delta v(t,x) \psi_R^2(x)\,\mathrm{d}x+\int_{\mathbb{R}^2}\ln(1+u(t,x)) \nabla v(t,x)\cdot \nabla \psi_R^2(x)\,\mathrm{d}x
\\
&+\int_{\mathbb{R}^2}u(t,x)\ln(1+u(t,x)) \nabla v(t,x)\cdot \nabla \psi_R^2(x)\,\mathrm{d}x.
\end{align*}
The third term of the right hand side further follows from the second equation in \eqref{eqn;KS} that
\begin{align*}
&-\int_{\mathbb{R}^2}u(t,x)\Delta v(t,x)\psi_R^2(x)\,\mathrm{d}x
\\
=\,&-\tau\int_{\mathbb{R}^2}\left[ \tau \partial_t v(t,x)-\Delta v(t,x)+\lambda v(t,x)\right]\partial_t v(t,x) \psi_R^2(x)\,\mathrm{d}x
\\
&-\lambda\int_{\mathbb{R}^2}u(t,x)v(t,x) \psi_R^2(x)\,\mathrm{d}x+\int_{\mathbb{R}^2} u^2(t,x) \psi_R^2(x)\,\mathrm{d}x
\\
\le\,&-\tau^2 \| \partial_t v \psi_R\|_2^2-\frac\tau2\frac{\mathrm{d}}{\mathrm{d}t}\int_{\mathbb{R}^2}|\nabla v(t,x)|^2\psi_R^2(x)\,\mathrm{d}x
\\
&-\tau\int_{\mathbb{R}^2}\partial_t v(t,x)\nabla v(t,x) \cdot \nabla \psi_R^2(x)\,\mathrm{d}x-\frac{\lambda\tau}{2}\frac{\mathrm{d}}{\mathrm{d}t}\int_{\mathbb{R}^2} v^2(t,x)\psi_R^2(x)\,\mathrm{d}x
+\int_{\mathbb{R}^2} u^2(t,x) \psi_R^2(x)\,\mathrm{d}x,
\end{align*}
and also the fifth term is estimated as
\begin{align*}
&\int_{\mathbb{R}^2}\ln(1+u(t,x)) \Delta v(t,x) \psi_R^2(x)\,\mathrm{d}x
\\
=\,&\int_{\mathbb{R}^2}\ln(1+u(t,x)) \left[\tau \partial_tv(t,x)+\lambda v(t,x)-u(t,x) \right] \psi_R^2(x)\,\mathrm{d}x
\\
\le\,&\tau\int_{\mathbb{R}^2}\sqrt{u(t,x)} |\partial_t v(t,x)|\psi_R^2(x)\,\mathrm{d}x+\lambda\int_{\mathbb{R}^2}\sqrt{u(t,x)} v(t,x) \psi_R^2(x)\,\mathrm{d}x
\\
\le\,&\frac{\tau^2}{4} \| \partial_t v \psi_R\|_2^2+C \| u_0\|_1+\lambda \|u_0\|_1^{\frac12}\sup_{t \ge 0}\|v(t)\|_2,
\end{align*}
where we use $\ln(1+x) \le \sqrt{x}$ for $x\ge0$.
Combining the above estimates implies that
\begin{align*}
&\frac{\mathrm{d}}{\mathrm{d}t}\mathcal{E}(t)+\int_{\mathbb{R}^2}\frac{|\nabla u(t,x)|^2}{1+u(t,x)} \psi_R^2(x)\,\mathrm{d}x+\frac{3\tau^2}4\|\partial_t v \psi_R\|_2^2
\\
\le\,&-\int_{\mathbb{R}^2}\ln(1+u(t,x)) \nabla u(t,x)\cdot \nabla \psi_R^2(x)\,\mathrm{d}x
-\int_{\mathbb{R}^2}\partial_t v(t,x)\nabla v(t,x) \cdot \nabla \psi_R^2(x)\,\mathrm{d}x
\\
&+\int_{\mathbb{R}^2} u^2(t,x) \psi_R^2(x)\,\mathrm{d}x
-\int_{\mathbb{R}^2}u(t,x) \nabla v(t,x)\cdot \nabla \psi_R^2(x)\,\mathrm{d}x
\\
&+\int_{\mathbb{R}^2}\ln(1+u(t,x)) \nabla v(t,x)\cdot \nabla \psi_R^2(x)\,\mathrm{d}x
\\
&+\int_{\mathbb{R}^2}u(t,x)\ln(1+u(t,x)) \nabla v(t,x)\cdot \nabla \psi_R^2(x)\,\mathrm{d}x+C,
\eqntag
\label{eqn;E}
\end{align*}
where
\begin{align*}
\mathcal{E}(t):=\int_{\mathbb{R}^2}(1+u(t,x))\ln(1+u(t,x))& \psi_R^2(x)\,\mathrm{d}x-\int_{\mathbb{R}^2}u(t,x)\psi_R^2(x)\,\mathrm{d}x
\\
&+\frac\tau2\int_{\mathbb{R}^2}|\nabla v(t,x)|^2\psi_R^2(x)\,\mathrm{d}x
+\frac{\tau\lambda}{2}\int_{\mathbb{R}^2} v^2(t,x)  \psi_R^2(x)\,\mathrm{d}x.
\end{align*}
We next show that all terms on the right hand side of \eqref{eqn;E} are estimated in terms of the $L^2$ norms of $u\psi_R$ and $|\nabla v| \chi_{ \supp\nabla\psi_R }$.
For the first term, since $\ln(1+x) \le \sqrt{x}$ for $x\ge0$ and 
\begin{align*}
\nabla \mathrm{H}(u(t,x)) = \ln(1+u(t,x))\nabla u(t,x),\qquad
\mathrm{H}(a):=\int_0^a \ln(1+s)\,\mathrm{d}s=(1+a)\ln(1+a)-a,
\end{align*}
 we obtain by Young's inequality for any $\varepsilon>0$,
\begin{align*}
&- \int_{\mathbb{R}^2}\ln(1+u(t,x)) \nabla u(t,x)\cdot \nabla \psi_R^2(x)\,\mathrm{d}x 
\\
=\,&\int_{\mathbb{R}^2}\mathrm{H}(u(t,x)) \Delta \psi_R^2(x)\,\mathrm{d}x
\\
\le\,&\frac{C}{R^2}\int_{\mathbb{R}^2} u(t,x)\sqrt{u(t,x)}\psi_R(x)\,\mathrm{d}x
\le\,\int_{\mathbb{R}^2}u^2(t,x)\psi_R^2(x)\,\mathrm{d}x+\frac{C}{R^4},
\end{align*}
and by applying Lemma \ref{lem;radial;heat-anul} as an annular region $A_{R/16}$,
\begin{align*}
\left|  \int_{\mathbb{R}^2}\partial_t v(t,x)\nabla v(t,x) \cdot \nabla \psi_R^2(x)\,\mathrm{d}x\right|
\le\,&\frac{\tau^2}{4}\|\partial_tv\psi_R\|_2^2+\frac{C}{R^2}\int_{\supp \nabla\psi_R}|\nabla v(t,x)|^2\,\mathrm{d}x
\\
\le\,&\frac{\tau^2}{4}\|\partial_tv\psi_R\|_2^2+\frac{C}{R^2}.
\end{align*}
As for the  fourth and fifth terms, by an argument similar to the preceding estimates
\begin{align*}
&\left| \int_{\mathbb{R}^2}u(t,x) \nabla v(t,x)\cdot \nabla \psi_R^2(x)\,\mathrm{d}x \right|
+\left|\int_{\mathbb{R}^2}\ln(1+u(t,x)) \nabla v(t,x)\cdot \nabla \psi_R^2(x)\,\mathrm{d}x \right|
\\
\le\,&\int_{\mathbb{R}^2}u^2(t,x)\psi_R^2(x)\,\mathrm{d}x+\frac{C}{R^2}\int_{\supp\nabla\psi_R}|\nabla v(t,x)|^2\,\mathrm{d}x+C\|u_0\|_1
\\
\le\,&\int_{\mathbb{R}^2}u^2(t,x)\psi_R^2(x)\,\mathrm{d}x+C \left(1+\frac{1}{R^2}\right).
\end{align*}
The last integral of \eqref{eqn;E} is estimated, by use of Lemma \ref{lem;radial2} with $\alpha=1$, as 
\begin{align*}
&\left| \int_{\mathbb{R}^2}u(t,x)\ln(1+u(t,x)) \nabla v(t,x)\cdot \nabla \psi_R^2(x)\,\mathrm{d}x \right|
\\
\le\,&\frac{C}{R} \left[ \frac{16\|u_0\|_1^\frac12}{\sqrt{2}\pi R}   \left( \int_{\mathbb{R}^2}\frac{|\nabla (u(t,x)\psi_R(x))|^2}{1+u(t,x)}\,\mathrm{d}x \right)^\frac12+1 \right]\int_{\supp\nabla \psi_R}\sqrt{u(t,x)}|\nabla v(t,x)|\,\mathrm{d}x
\\
\le\,&\frac{C\|u_0\|_1^{\frac12}}{R} \left[ \frac{16\|u_0\|_1^\frac12}{\sqrt{2}\pi R} \left( \int_{\mathbb{R}^2}\frac{|\nabla (u(t,x)\psi_R(x))|^2}{1+u(t,x)}\,\mathrm{d}x \right)^\frac12+1 \right]
\left( \int_{\supp\nabla \psi_R}|\nabla v(t,x)|^2\,\mathrm{d}x\right)^{\frac12}
\\
\le\,&\varepsilon\int_{\mathbb{R}^2}\frac{|\nabla u(t,x)|^2}{1+u(t,x)}\psi_R^2(x)\,\mathrm{d}x+C(\varepsilon)\left(1+\frac{1}{R^4}+\frac{1}{R^2}\right).
\end{align*}
Hence, it follows from \eqref{eqn;E} together with the above estimates that for $R>1$
\begin{align*}
\frac{\mathrm{d}}{\mathrm{d}t}\mathcal{E}(t)+\mathcal{E}(t)+(1-2\varepsilon)&\int_{\mathbb{R}^2}\frac{|\nabla u(t,x)|^2}{1+u(t,x)} \psi_R^2(x)\,\mathrm{d}x+\frac{\tau^2}2\|\partial_t v \psi_R\|_2^2
\\
\le\,&4
\int_{\mathbb{R}^2} u^2(t,x) \psi_R^2(x)\,\mathrm{d}x
+\frac{\tau}{2}
\int_{\mathbb{R}^2}|\nabla v(t,x)|^2\psi_R^2(x)\,\mathrm{d}x
+C(\varepsilon),
\end{align*}
where $(1+x)\ln(1+x) -x\le x^2$ for $x\ge0$ and $\sup_{t \ge 0}\|v(t)\|_2<\infty$ for $\lambda>0$.
Moreover, Lemma \ref{lem;radial;heat-anul} along with the $L^2$ bound on $v$ imply that
\begin{align*}
\frac\tau2\int_{\mathbb{R}^2}|\nabla v(t,x)|^2\psi_R^2(x)\,\mathrm{d}x=\,&-\frac\tau2\int_{\mathbb{R}^2}v(t,x)\Delta v(t,x) \psi_R^2(x)\,\mathrm{d}x-\frac\tau2\int_{\mathbb{R}^2}v(t,x)\nabla v(t,x)\cdot \nabla\psi_R^2(x)\,\mathrm{d}x
\\
=\,&-\frac{\tau^2}{2}\int_{\mathbb{R}^2}v(t,x)\partial_t v(t,x)\psi_R^2(x)\,\mathrm{d}x-\frac{\tau\lambda}{2}\int_{\mathbb{R}^2}v^2(t,x)\psi_R^2(x)\,\mathrm{d}x
\\
&+\frac{\tau}{2}\int_{\mathbb{R}^2}v(t,x)u(t,x)\psi_R^2(x)\,\mathrm{d}x
-\frac\tau2\int_{\mathbb{R}^2}v(t,x)\nabla v(t,x)\cdot \nabla\psi_R^2(x)\,\mathrm{d}x
\\
\le\,&\frac{\tau^2}{4} \| \partial_t v \psi_R\|_2^2+\int_{\mathbb{R}^2}u^2(t,x) \psi_R^2(x)\,\mathrm{d}x+ C.
\end{align*}
Hence, we deduce that
\begin{align*}
\frac{\mathrm{d}}{\mathrm{d}t}\mathcal{E}(t)+\mathcal{E}(t)+(1-2\varepsilon)\int_{\mathbb{R}^2}\frac{|\nabla u(t,x)|^2}{1+u(t,x)} \psi_R^2(x)\,\mathrm{d}x&+\frac{\tau^2}4\|\partial_t v \psi_R\|_2^2
\\
\le\,&5
\int_{\mathbb{R}^2} u^2(t,x) \psi_R^2(x)\,\mathrm{d}x
+C(\varepsilon).
\end{align*}
Finally, we apply Lemma \ref{lem;radial2} with $\alpha=1$ to obtain
\begin{align*}
5
\int_{\mathbb{R}^2} u^2(t,x) \psi_R^2(x)\,\mathrm{d}x
\le\,&5\|u_0\|_1\left[ \frac{16}{\sqrt{2}\pi R}  \|u_0\|_1^\frac12 \left( \int_{\mathbb{R}^2}\frac{|\nabla (u(t,x) \psi_R(x))|^2}{1+u(t,x)}\,\mathrm{d}x \right)^\frac12+1  \right]
\\
\le\,&\varepsilon \int_{\mathbb{R}^2}\frac{|\nabla u(t,x) |^2}{1+u(t,x)}\psi_R^2(x)\,\mathrm{d}x + C(\varepsilon,R),
\end{align*}
which leads to
\begin{align*}
\frac{\mathrm{d}}{\mathrm{d}t}\mathcal{E}(t)+\mathcal{E}(t)+\frac12\int_{\mathbb{R}^2}\frac{|\nabla u(t,x)|^2}{1+u(t,x)} \psi_R^2(x)\,\mathrm{d}x+\frac{\tau^2}4\|\partial_t v \psi_R\|_2^2
\le\,&C
\end{align*}
by choosing $\varepsilon=1/6$.
As a consequence, Gronwall's inequality implies that
\begin{align*}
\sup_{t \ge t_0}\left[  \int_{\mathbb{R}^2}(1+u(t,x))\ln(1+u(t,x)) \psi_R^2(x)\,\mathrm{d}x+\int_{\mathbb{R}^2}|\nabla v(t,x)|^2\psi_R^2(x)\,\mathrm{d}x\right] <\infty.
\end{align*}
Accordingly,
\begin{align*}
\sup_{t \ge t_0}\int_t^{t+1} \left[ \int_{\mathbb{R}^2}\frac{ |\nabla u(s,x)|^2}{1+u(s ,x)} \psi_R^2(x)\,\mathrm{d}x+ \|\partial_tv(s) \psi_R\|_2^2 \right]\,\mathrm{d}s<\infty,
\end{align*}
which concludes the proof.
\end{pr}

\begin{lemma}\label{lem;L2L2-u}
Under the assumption as in Lemma \ref{lem;out-log}, 
\begin{align*}
\sup_{t \ge t_0}\int_t^{t+1}\left(  \int_{\mathbb{R}^2}u^2(s,x)\psi_{2R}^2(x)\,\mathrm{d}x+ \int_{|x|\ge \frac R8} \left[ u^2(s,x)+ |\partial_t v(s,x)|^2\right]\,\mathrm{d}x\right)\mathrm{d}s<\infty.
\end{align*}
\end{lemma}
\begin{pr}{Lemma \ref{lem;L2L2-u}}
Noting that  $\psi_R\equiv1$ on $|x|\ge R/8$,
we observe that
the claim is readily obtained from Lemma~\ref{lem;out-log} with the help of Lemma \ref{lem;radial2}.\end{pr}

\begin{lemma}\label{lem;L2-u}
Suppose Lemma \ref{lem;L2L2-u} holds.
Then,
\begin{align*}
\sup_{t \ge t_0} \left[ \| u(t)\psi_{2R}\|_2+ \int_t^{t+1} \int_{\mathbb{R}^2} |\nabla u(s,x)|^2\psi_{2R}^2(x)\,\mathrm{d}x\mathrm{d}s\right] <\infty.
\end{align*}
\end{lemma}
\begin{pr}{Lemma \ref{lem;L2-u}}
We begin with a standard energy estimate.
Using the first equation in~\eqref{eqn;KS}, we obtain 
\begin{align*}
&\frac12\frac{\mathrm{d}}{\mathrm{d}t} \int_{\mathbb{R}^2} u^2(t,x)\psi_{2R}^2(x)\,\mathrm{d}x\\
=\,&-\int_{\mathbb{R}^2}|\nabla u(t,x)|^2\psi_{2R}^2(x)\,\mathrm{d}x
-\int_{\mathbb{R}^2} u(t,x)\nabla u(t,x)\cdot \nabla \psi_{2R}^2(x)\,\mathrm{d}x
\\
&+\int_{\mathbb{R}^2}u(t,x) \nabla u(t,x)\cdot \nabla v(t,x) \psi_{2R}^2(x)\,\mathrm{d}x+\int_{\mathbb{R}^2}u^2(t,x)\nabla v(t,x)\cdot \nabla \psi_{2R}^2(x)\,\mathrm{d}x
\\
=\,&-\int_{\mathbb{R}^2}|\nabla u(t,x)|^2\psi_{2R}^2(x)\,\mathrm{d}x+\frac12\int_{\mathbb{R}^2} u^2(t,x)\Delta \psi_{2R}^2(x)\,\mathrm{d}x
\\
&-\frac12\int_{\mathbb{R}^2}u^2(t,x)\Delta v (t,x) \psi_{2R}^2(x)\,\mathrm{d}x+\frac12\int_{\mathbb{R}^2}u^2(t,x)\nabla v(t,x)\cdot \nabla \psi_{2R}^2(x)\,\mathrm{d}x
\\
\le\,&-\int_{\mathbb{R}^2}|\nabla u(t,x)|^2\psi_{2R}^2(x)\,\mathrm{d}x+\frac12\int_{\mathbb{R}^2} u^2(t,x)\Delta \psi_{2R}^2(x)\,\mathrm{d}x
\\
&-\frac\tau2\int_{\mathbb{R}^2}u^2(t,x)\partial_tv (t,x) \psi_{2R}^2(x)\,\mathrm{d}x+\frac{1}2\int_{\mathbb{R}^2}u^{3}(t,x) \psi_{2R}^2(x)\,\mathrm{d}x
\\
&+\frac12\int_{\mathbb{R}^2}u^2(t,x)\nabla v(t,x)\cdot \nabla \psi_{2R}^2(x)\,\mathrm{d}x.
\end{align*}
Lemma \ref{lem;radial-lemma} with $(p,q)=(2,2)$ and Young's inequality imply that
\begin{align*}
\left| \frac12\int_{\mathbb{R}^2} u^2(t,x)\Delta \psi_{2R}^2(x)\,\mathrm{d}x  \right|
\le\,&\frac{C}{R^2}\int_{\supp \psi_{2R}}u^2(t,x)\psi_{2R}\,\mathrm{d}x
\\
\le\,&\frac{C \|u_0\|_1}{R^2} \left( \frac{8}{\pi R} \right)^{\frac12} \| u\psi_{2R} \|_2^{\frac12} \|\nabla(u\psi_{2R})\|_2^{\frac12}
\\
\le\,&\frac18 \| |\nabla u|\psi_{2R}\|_2^2+C(R)\|u\psi_{2R}\|_2^2+\frac{C}{R^5}.
\end{align*}
Similarly,
\begin{align*}
\left|  \frac12\int_{\mathbb{R}^2}u^2(t,x)\partial_tv (t,x) \psi_{2R}^2(x)\,\mathrm{d}x  \right|
\le\,&\frac12 \left( \frac{8}{\pi R} \right)^{\frac12} \| u\psi_{2R} \|_2^{\frac12} \|\nabla(u\psi_{2R})\|_2^{\frac12}\| u\psi_{2R}\|_2 \|\partial_t v\|_{ L^2(|x|>\frac R8)}
\\
\le\,&\frac18\| |\nabla u|\psi_{2R}\|_2^2+C(R)\|u\psi_{2R}\|_2^2+\frac{C}{R}\|u\psi_{2R}\|_2^2\|\partial_t v\|_{  L^2(|x|>\frac R8)}^2,
\end{align*}
as well as it follows from $\supp\nabla \psi_{2R}\subset \{ R/8 \le |x| \le R/4 \}$ and Lemma \ref{lem;radial;heat-anul} that
\begin{align*}
\left|  \frac12\int_{\mathbb{R}^2}u^2(t,x)\nabla v(t,x)\cdot \nabla \psi_{2R}^2(x)\,\mathrm{d}x \right|
\le\,&\frac{C}{R}\int_{\frac{R}{8}\le |x|\le\frac R4}u^2(t,x) |\nabla v(t,x)| \psi_{2R}(x)\,\mathrm{d}x
\\
\le\,&\frac{C}{R} \left( \frac{8}{\pi R} \right)^{\frac12} \| u\psi_{2R} \|_2^{\frac12} \|\nabla(u\psi_{2R})\|_2^{\frac12}\|u\|_{ L^2( |x|\ge \frac R8 ) } \| \nabla v\|_{ L^2( \frac R{8}\le|x|\le\frac R4 ) }
\\
\le\,&\frac18\| |\nabla u| \psi_{2R}\|_2^2+C(R) \| u\psi_{2R}\|_2^2+\frac{C}{R^3} \| u\|_{  L^2(|x|\ge\frac R8) }^2.
\end{align*}
Finally, since we estimate by Lemma \ref{lem;radial-lemma} and Young's inequality
\begin{align*}
\left|  \frac{1}2\int_{\mathbb{R}^2}u^{3}(t,x) \psi_{2R}^2(x)\,\mathrm{d}x \right|
\le\,&\frac{C\|u_0\|_1}{R}   \| u\psi_{2R} \|_2 \|\nabla(u\psi_{2R})\|_2
\\
\le\,& \frac18\| |\nabla u| \psi_{2R}\|_2^2+ C(R) \| u \psi_{2R}\|_2^2,
\end{align*}
combining above computations arrives at
\begin{align*}
\frac12\frac{\mathrm{d}}{\mathrm{d}t}\|u\psi_{2R}\|_2^2+\frac12\| |\nabla u| \psi_{2R}\|_2^2
\le\,&C(R) \|u\psi_{2R}\|_2^2\left( 1+\|\partial_t v\|_{ L^2(|x|\ge \frac R8) }^2 \right)
+C(R) \left(\| u\|_{  L^2(|x|\ge\frac R8) }^2+1 \right).
\end{align*}
Therefore, uniform Gronwall's inequality (cf.~\cite[Lemma $\mathrm{III}$.1.1]{Te}) with the help of Lemma \ref{lem;L2L2-u} gives the desired estimates,
which completes the proof.
\end{pr}

\begin{lemma}\label{lem;lp-u}
Suppose that Lemma \ref{lem;L2-u} holds. Then, for any $p\in(2,\infty)$
\begin{align*}
\sup_{t \ge t_0}\| u\psi_{4R} \|_p+\sup_{t \ge t_0}\int_t^{t+1} \int_{\mathbb{R}^2} | \nabla u^\frac p2(s,x)|^2\psi_{4R}^p(x)\,\mathrm{d}x\mathrm{d}s<\infty.
\end{align*}
\end{lemma}
\begin{pr}{Lemma \ref{lem;lp-u}}
The proof follows a similar argument to that of Lemma~\ref{lem;L2-u}.
Let $p\in (2,\infty)$ and $R>1$. 
Since it follows from the second equation in \eqref{eqn;KS} that
\begin{align*}
\tau\partial_t(v\psi_{2R})=\Delta(v\psi_{2R})-\lambda(v\psi_{2R})-2\nabla v\cdot \nabla \psi_{2R} - v \Delta \psi_{2R} + u\psi_{2R},
\end{align*}
the integral formulation is given as
\begin{align*}
(v\psi_{2R})(t)=\,&e^{(t-t_0) \tau^{-1}(\Delta-\lambda)}(v(t_0)\psi_{2R})
\\&+\frac1\tau\int_{t_0}^t e^{(t-s) \tau^{-1}(\Delta-\lambda)}\left[ -2\nabla v\cdot \nabla \psi_{2R} - v \Delta \psi_{2R} + u\psi_{2R}  \right](s)\,\mathrm{d}s,
\eqntag
\label{eqn;vpsi}
\end{align*}
from which we derive from Lemma \ref{lem;L2-u} along with Lemma \ref{lem;radial;heat-anul} that
\begin{align*}
\|\nabla (v\psi_{2R})(t)\|_p\le\,&\|\nabla (v(t_0)\psi_{2R})\|_p
\\
&+C\int_{t_0}^t e^{-\lambda\tau^{-1}(t-s)}(t-s)^{-1+\frac1p}\left[ \| |\nabla v(s)| |\nabla \psi_{2R}|\|_2+\|v(s) \Delta \psi_{2R}\|_2+\|u(s)\psi_{2R}\|_2  \right]\,\mathrm{d}s
\\
\le\,&\|\nabla (v(t_0)\psi_{2R})\|_p+C(\lambda,p) \Gamma\left(\frac1p\right),\qquad t\ge t_0.
\end{align*}
Hence, by use of  the integrability $\sup_{t \ge 0}\|v(t)\|_p<\infty$ for $\lambda>0$, we have
\begin{align*}
\sup_{t \ge t_0} \||\nabla v(t)|\psi_{2R}\|_p <\infty,
\end{align*}
and accordingly
\begin{align*}
\sup_{t \ge t_0}\| \nabla v(t)\|_{ L^p(  |x| \ge \frac R4) }<\infty.
\eqntag
\label{eqn;nabla-v-p}
\end{align*}
We here recall that $\psi_{4R}\equiv1$ on $|x|\ge R/2$ and $\supp \psi_{4R}\subset \{x\in\mathbb{R}^2; |x|\ge R/4\}$.
Then, it follows from integrations by parts that
\begin{align*}
&\frac1p\frac{\mathrm{d}}{\mathrm{d}t} \int_{\mathbb{R}^2} u^p(t,x)\psi_{4R}^p(x)\,\mathrm{d}x\\
=\,&-\frac{4(p-1)}{p^2}\int_{\mathbb{R}^2}|\nabla u^\frac p2(t,x)|^2\psi_{4R}^p(x)\,\mathrm{d}x
-\int_{\mathbb{R}^2} u^{p-1}(t,x)\nabla u(t,x)\cdot \nabla \psi_{4R}^p(x)\,\mathrm{d}x
\\
&+\int_{\mathbb{R}^2}u(t,x) \nabla u^{p-1}(t,x)\cdot \nabla v(t,x) \psi_{4R}^p(x)\,\mathrm{d}x+\int_{\mathbb{R}^2}u^p(t,x)\nabla v(t,x)\cdot \nabla \psi_{4R}^p(x)\,\mathrm{d}x
\\
=\,&-\frac{4(p-1)}{p^2}\int_{\mathbb{R}^2}|\nabla u^\frac p2(t,x)|^2\psi_{4R}^p(x)\,\mathrm{d}x+\frac1p\int_{\mathbb{R}^2} u^p(t,x)\Delta \psi_{4R}^p(x)\,\mathrm{d}x
\\
&+\frac{2(p-1)}p\int_{\mathbb{R}^2}u^\frac p2(t,x)\nabla u^\frac p2(t,x)\cdot \nabla v (t,x) \psi_{4R}^p(x)\,\mathrm{d}x+\int_{\mathbb{R}^2}u^p(t,x)\nabla v(t,x)\cdot \nabla \psi_{4R}^p(x)\,\mathrm{d}x.
\end{align*}
Using Lemma \ref{lem;radial-lemma}, the estimate
$|\nabla(\psi_{4R}^{p/2})|^2 \le C_p R^{-2}\psi_{4R}^{p-1}$,
and Young's inequality, we obtain
\begin{align*}
\int_{ \mathbb{R}^2} u^p(t,x) \psi_{4R}^{p-1}(x)\,\mathrm{d}x\le\,&
\|u_0\|_1\left[  \left(\frac{2(p+1)}{\pi p R}\right)^{\frac{2}{p+1}} \| u\psi_{4R}\|_1^{\frac{1}{p+1}} \|\nabla(u\psi_{4R})^{\frac p2}\|_2^{\frac{2}{p+1}} \right]^{p-1}
\\
\le\,&\frac{C}{R^{\frac{2(p-1)}{p+1}}}  \left[  \int_{ \mathbb R^2} |\nabla u^\frac p2(t,x)|^2\psi_{4R}^p\,\mathrm{d}x +\frac{1}{R^2}\int_{ \mathbb R^2} u^p(t,x) \psi_{4R}^{p-1}(x)\,\mathrm{d}x \right]^{\frac{p-1}{p+1}}
\\
\le\,&\frac{C}{R^{\frac{2(p-1)}{p+1}}}  \left(  \int_{ \mathbb R^2} |\nabla u^\frac p2(t,x)|^2\psi_{4R}^p\,\mathrm{d}x   \right)^{\frac{p-1}{p+1}}+\frac{C}{R^{\frac{4(p-1)}{p+1}}} \left(  \int_{ \mathbb R^2} u^p(t,x) \psi_{4R}^{p-1}(x)\,\mathrm{d}x \right)^{\frac{p-1}{p+1}}
\\
\le\,&\frac{C}{R^{\frac{2(p-1)}{p+1}}}  \left(  \int_{ \mathbb R^2} |\nabla u^\frac p2(t,x)|^2\psi_{4R}^p\,\mathrm{d}x   \right)^{\frac{p-1}{p+1}}+\frac12  \int_{ \mathbb R^2} u^p(t,x) \psi_{4R}^{p-1}(x)\,\mathrm{d}x+C(R,p).
\end{align*}
Absorbing the second term on the right-hand side into the
left-hand side and applying Young's inequality once more,
we deduce that, for any $\varepsilon>0$,
\begin{align*}
\frac{C}{R^2} \int_{ \mathbb R^2} u^p(t,x) \psi_{4R}^{p-1}(x)\,\mathrm{d}x\le\,&\frac{C}{R^{2+\frac{2(p-1)}{p+1}}}  \left(  \int_{ \mathbb R^2} |\nabla u^\frac p2(t,x)|^2\psi_{4R}^p\,\mathrm{d}x   \right)^{\frac{p-1}{p+1}}+C(R,p)
 \\
 \le\,&\varepsilon\int_{ \mathbb R^2} |\nabla u^\frac p2(t,x)|^2\psi_{4R}^p\,\mathrm{d}x+C(\varepsilon,R,p).
 \eqntag
 \label{eqn;lp-u}
\end{align*}
Hence, as for the second term, 
\begin{align*}
\left|\frac1p\int_{\mathbb{R}^2} u^p(t,x)\Delta \psi_{4R}^p(x)\,\mathrm{d}x \right|
\le\,&\frac{C}{R^2}\int_{\mathbb{R}^2 } u^p(t,x) \psi_{4R}^{p-1}(x)\,\mathrm{d}x
\\
\le\,&\varepsilon\| |\nabla u^{\frac p2}| \psi_{4R}^\frac p2\|_2^2
+C(\varepsilon,R).
\end{align*}
Next, applying Young's and Gagliardo--Nirenberg's inequalities together with \eqref{eqn;nabla-v-p} gives that
\begin{align*}
&\left| \frac{2(p-1)}p\int_{\mathbb{R}^2}u^\frac p2(t,x)\nabla u^\frac p2(t,x)\cdot \nabla v (t,x) \psi_{4R}^p(x)\,\mathrm{d}x \right|
\\
\le\,&\frac{p-1}{2p^2}\| |\nabla u^{\frac p2}|\psi_{4R}^\frac p2\|_2^2
+C(R,p)\int_{ \mathbb R^2}u^p(t,x)|\nabla v(t,x)|^2\psi_{4R}^p(x)\,\mathrm{d}x
\\
\le\,&\frac{p-1}{2p^2}\| |\nabla u^{\frac p2}|\psi_{4R}^\frac p2\|_2^2+C(R,p)\|(u\psi_{4R})^\frac p2\|_4^2\|\nabla v\|_{L^4(\supp\psi_{4R})}^2
\\
\le\,&\frac{p-1}{2p^2}\| |\nabla u^{\frac p2}|\psi_{4R}^\frac p2\|_2^2+C(R,p)\|(u\psi_{4R})^\frac p2\|_2\|\nabla(u\psi_{4R})^\frac p2\|_2
\\
\le\,&\left(  \frac{p-1}{p^2}+\varepsilon\right)\| |\nabla u^{\frac p2}|\psi_{4R}^\frac p2\|_2^2+C(R,p)\|u\psi_{4R}\|_p^p+ C(\varepsilon,R,p).
\end{align*}
The last term is similarly estimated from Lemma \ref{lem;L2-u} and Lemma \ref{lem;radial;heat-anul} as
\begin{align*}
&\left|  \int_{\mathbb{R}^2}u^p(t,x)\nabla v(t,x)\cdot \nabla \psi_{4R}^p(x)\,\mathrm{d}x\right|
\\
\le\,&\frac{C}{R}\int_{ \supp\nabla \psi_{4R}} u^p(t,x)|\nabla v(t,x)| \psi_{4R}^{p-1}\,\mathrm{d}x
\\
\le\,&\frac{C}{R} \|u\|_{L^2(\supp\nabla\psi_{4R})} \|\nabla v\|_{L^2(\supp\nabla\psi_{4R})}\left[  \left(\frac{2(p+1)}{\pi p R}\right)^{\frac{2}{p+1}} \| u\psi_{4R}\|_1^{\frac{1}{p+1}} \|\nabla(u\psi_{4R})^{\frac p2}\|_2^{\frac{2}{p+1}} \right]^{p-1}
\\
\le\,&\varepsilon\| |\nabla u^{\frac p2}|\psi_{4R}^\frac p2\|_2^2+ C(\varepsilon,R,p),
\end{align*}
where $\supp\nabla \psi_{4R}\subset \{x\in\mathbb{R}^2; R/4\le|x|\le R/2 \}$.
Hence, we end up with 
\begin{align*}
\frac1p\frac{\mathrm{d}}{\mathrm{d}t} \|u\psi_{4R}\|_p^p+\frac1p\|u\psi_{4R}\|_p^p
+\left(\frac{3(p-1)}{p^2}-3\varepsilon\right)\| |\nabla u^{\frac p2}|\psi_{4R}^\frac p2\|_2^2
\le\,&C(R)\|u\psi_{4R}\|_p^p+C(\varepsilon,R)
\end{align*}
for $t\ge t_0$.
Letting $3\varepsilon=(p-1)/p^2$, we have
\begin{align*}
\frac1p\frac{\mathrm{d}}{\mathrm{d}t} \|u\psi_{4R}\|_p^p+\frac1p\|u\psi_{4R}\|_p^p
+\frac{2(p-1)}{p^2}\| |\nabla u^{\frac p2}|\psi_{4R}^\frac p2\|_2^2
\le\,&C(R)\|u\psi_{4R}\|_p^p+C(R,p).
\end{align*}
According to \cite[Lemma 2.4]{Su06} with $(q_1,q_2,r)=(1,p,p)$ and $m=1$,
\begin{align*}
\|u\psi_{4R}\|_p^p\le\,&C \|u\psi_{4R}\|_1^{p(1-\theta)}\|\nabla(u\psi_{4R})^\frac p2\|_2^{2\theta},
\qquad
\theta=1-\frac1p,
\end{align*}
from which it follows from Young's inequality, mass conservation and \eqref{eqn;lp-u} that eventually
\begin{align*}
C(R)\|u\psi_{4R}\|_p^p\le\,&C\|\nabla(u\psi_{4R})^\frac p2\|_2^{2\theta}
\le\,\frac{p-1}{p^2}\| |\nabla u^\frac p2|\psi_{4R}^{\frac p2}\|_2^2+C(R,p).
\end{align*}
Therefore,
we obtain by Gronwall's inequality,
\begin{align*}
\sup_{t \ge t_0}\|u(t)\psi_{4R}\|_p+\sup_{t \ge t_0}\int_t^{t+1} \| |\nabla u^\frac p2(s)| \psi_{4R}^\frac p2\|_2^2\,\mathrm{d}s<\infty,
\qquad
p>2,
\end{align*}
which concludes the proof.
\end{pr}

\vspace{3mm}
We are in a position to see the $L^\infty$ outer estimates for $u$.
\begin{lemma}\label{lem;infty-u}
Suppose that Lemma \ref{lem;lp-u} holds. Then,
\begin{align*}
\sup_{t \ge t_0}\left( \|u(t) \psi_{8R}\|_\infty+\| \nabla (v(t) \psi_{4R})\|_\infty\right)<\infty.
\end{align*}
Consequently,
\begin{align*}
\sup_{t \ge t_0}\| u(t)\|_{ L^\infty(|x|\ge R) }<\infty.
\end{align*}
\end{lemma}
\begin{pr}{Lemma \ref{lem;infty-u}}
First of all, notice that
by virtue of Lemma \ref{lem;lp-u} and \eqref{eqn;nabla-v-p}, we have the uniform-in-time $L^p$ outer estimates for $u\psi_{4R}$ and $\nabla v$ that
\begin{align*}
\sup_{t \ge t_0}\|u(t)\psi_{4R}\|_p<\infty,\qquad
\sup_{t \ge t_0}\| \nabla v(t)\|_{ L^p(  |x| \ge \frac R4) }<\infty,\qquad p\in(2,\infty),
\eqntag
\label{eqn;uv-lp}
\end{align*}
where $\psi_{4R}$ is the cut-off defined in \eqref{eqn;cutoff}.
Recalling the formulation \eqref{eqn;vpsi} as $\psi_{4R}$,
\begin{align*}
(v\psi_{4R})(t)=\,&e^{(t-t_0) \tau^{-1}(\Delta-\lambda)}(v(t_0)\psi_{4R})
\\&+\frac1\tau\int_{t_0}^t e^{(t-s) \tau^{-1}(\Delta-\lambda)}\left[ -2\nabla v\cdot \nabla \psi_{4R} - v \Delta \psi_{4R} + u\psi_{4R}  \right](s)\,\mathrm{d}s,
\end{align*}
we obtain
\begin{align*}
\| \nabla(v\psi_{4R})(t)\|_\infty\le\,&\|\nabla (v(t_0)\psi_{4R})\|_\infty
\\
&+C\int_{t_0}^t e^{-(t-s)\frac{\lambda}{\tau}}(t-s)^{-\frac1p-\frac12} \left[ \| |\nabla v(s)| |\nabla \psi_{4R}|\|_p+ \|v(s)\Delta\psi_{4R}\|_p+\|u(s)\psi_{4R}\|_p \right]\,\mathrm{d}s
\\
\le\,&\|\nabla (v(t_0)\psi_{4R})\|_\infty+C\Gamma\left(\frac12-\frac1p\right)
\end{align*}
for $t\ge t_0$ and $p>2$.
Hence,
\begin{align*}
\sup_{t \ge t_0} \| \nabla(v\psi_{4R})(t)\|_\infty<\infty.
\end{align*}
Next, by use of the first equation in \eqref{eqn;KS}
\begin{align*}
\partial_t(u\psi_{8R} )=\,&(\Delta-1)(u\psi_{8R})+u\psi_{8R}+u(\Delta\psi_{8R}+\nabla v\cdot \nabla \psi_{8R})-2\nabla\cdot(u\nabla \psi_{8R})-\nabla\cdot(u\psi_{8R}\nabla v)
\\
=\,&(\Delta-1)(u\psi_{8R})+F(u,v,\psi_{8R}),
\eqntag
\label{eqn;1st-exterior-eq}
\end{align*}
and we have by use of \eqref{eqn;1st-exterior-eq}
\begin{align*}
(u\psi_{8R})(t)=\,&e^{(t-t_0)(\Delta-1)}((u\psi_{8R})(t_0))+\int_{t_0}^t  e^{(t-s) (\Delta-1)} F(u(s),v(s),\psi_{8R})\,\mathrm{d} s.
\end{align*}
Since \eqref{eqn;uv-lp} and
\begin{align*}
\sup_{t \ge t_0}\|u(t)\|_{L^p(|x|\ge \frac R2)} \le\, \sup_{t \ge t_0}\|u(t)\psi_{4R}\|_p<\infty,
\end{align*}
it follows from  the heat semigroup estimates and $\supp\nabla\psi_{8R}\subset\{x\in\mathbb{R}^2; R/2\le|x|\le R\}$ that
for~$p>2$
\begin{align*}
&\left \|\int_{t_0}^t e^{(t-s) (\Delta-1)}F(u(s),v(s),\psi_{8R})\,\mathrm{d} s\right\|_\infty
\\
\le\,&C\int_{t_0}^t e^{-(t-s)} (t-s)^{-\frac1p} \left(\|u(s)\psi_{8R}\|_p+\|u(s) (\Delta\psi_{8R}+\nabla v(s)\cdot\nabla\psi_{8R})\|_p\right)\,\mathrm{d}s
\\
&+C\int_{t_0}^t e^{-(t-s)} (t-s)^{-\frac1p-\frac12} \left( \|u(s)|\nabla\psi_{8R}|\|_p + \|u(s)\psi_{8R}|\nabla v|\|_p\right)\,\mathrm{d}s
\\
\le\,&C(R,p)\Gamma\left(1-\frac1p\right)+C(R,p) \Gamma\left(\frac12-\frac1p\right),
\end{align*}
so that
\begin{align*}
\sup_{t \ge t_0}\|(u\psi_{8R})(t)\|_\infty
\le\,&\|(u\psi_{8R})(t_0)\|_\infty+C(R,p).
\end{align*}
This finally implies that
\begin{align*}
\sup_{t \ge t_0}\| u(t)\|_{ L^\infty(|x|\ge R) }\le\,&\sup_{t \ge t_0}\|(u\psi_{8R})(t)\|_\infty<\infty,
\end{align*}
which ends the proof.
\end{pr}


\vspace{5mm}
\section{Uniform-in-time interior estimates}\label{sect;3}
In this section, we derive interior estimates for radially symmetric solutions to \eqref{eqn;KS} under the optimal mass constraint \eqref{eqn;mass-constraint}.
As discussed in the introduction, the classical free energy argument on the whole space does not directly yield uniform-in-time boundedness, since the available estimate for the negative part of the entropy may grow  in time.
To overcome this difficulty, we consider a localized free energy functional on a ball of radius $R$, chosen so that the exterior estimates established in Section \ref{sect;2} apply near its boundary.
The key ingredient is Lemma \ref{lem;log}, which controls the interior entropy in terms of the local energy dissipation.
Combined with the exterior estimates, this yields uniform-in-time bounds in the interior region.
It is worth pointing out that a radius $R$  in all statements below can be taken arbitrarily, and also all constants do not depend on the time $t$.

\subsection{Localized energy inequality for uniform-in-time estimates}
Fix $R>0$.
Let $\mathcal{F}_R$ be a localized free energy defined by
\begin{align*}
\mathcal{F}_R(t):=\int_{\mathbb{R}^2}u(t,x)(\ln u(t,x)-1)\varphi_R^2(x)\,\mathrm{d}x-&\int_{\mathbb{R}^2}u(t,x)v(t,x)\varphi_R^2(x)\,\mathrm{d}x
\\
&+\frac12\int_{\mathbb{R}^2}|\nabla v(t,x)|^2\varphi_R^2(x)\,\mathrm{d}x+\frac{\lambda}{2}\int_{\mathbb{R}^2} v^2(t,x)\varphi_R^2(x)\,\mathrm{d}x,
\end{align*}
where $\varphi_R:=1-\psi_{16R}$.
We first mention our main goal in this subsection.
\begin{proposition}\label{prop1}
Let $(u,v)$ be the radially symmetric solution to \eqref{eqn;KS}.
Assume  \eqref{eqn;mass-constraint} and that for any $R>0$ and $t_0>0$,  there exists $C(R)=C(R,t_0)>0$ such that 
\begin{align*}
\sup_{t \ge t_0}\| \sqrt{u(t) } (\ln u(t)-v)\|_{L^2(R< |x|<2R )}+ \sup_{t \ge t_0}u(t,R)\le\, C(R).
\eqntag
\label{eqn;assumption-R}
\end{align*}
Then, there are $\alpha$, $\beta>0$ depending on $R$ and $t_0$ such that
\begin{align*}
\frac{\mathrm{d}}{\mathrm{d}t} \mathcal{F}_R(t)+\alpha \mathcal{F}_R(t)+ \beta \mathcal{D}_R(t)\le\, C(R),
\eqntag
\label{eqn;energy-inequality}
\end{align*} 
where 
\begin{align*}
\mathcal{D}_R(t):=\int_{\mathbb{R}^2 }u(t,x) |\nabla(\ln u(t,x)  -v(t,x) )|^2\varphi_R^2(x)\,\mathrm{d}x+\tau \int_{ \mathbb{R}^2}|\partial_tv(t,x) |^2\varphi_R^2(x) \,\mathrm{d}x.
\end{align*}
\end{proposition}
\vspace{3mm}
\begin{remark}
Once the energy inequality \eqref{eqn;energy-inequality}
is established, Gronwall's inequality yields a uniform-in-time
upper bound for $\mathcal{F}_R$.
The existence of a radius $R>1$ satisfying
\eqref{eqn;assumption-R} follows from the exterior estimates
established in Section \ref{sect;2}. Indeed, Lemma \ref{lem;infty-u} implies that
\begin{align*}
\sup_{t \ge t_0}u(t,R) \le\,\sup_{t \ge t_0} \|u(t)\|_{L^\infty(|x|\ge R)} < \infty.
\end{align*}
Deriving \eqref{eqn;energy-inequality} is nontrivial, since
the dissipative term
$\|\sqrt{u}\nabla(\ln u-v)\|_2^2$
associated with the first equation in \eqref{eqn;KS} does not directly control
the standard Fisher information
$\|\sqrt{u}\nabla\ln u\|_2^2$ due to lack of regularity of solutions.
Thus, additional arguments are needed to recover
the regularity estimates required for boundedness.
Differential energy inequalities of this type play a crucial role
in establishing boundedness of solutions, see for instance
\cite{Bi-Ci-Fu-Se} for a one-dimensional critical quasilinear fully parabolic
Keller--Segel system on an interval.
However, the argument in \cite{Bi-Ci-Fu-Se} relies heavily
on identities and estimates specific to one space dimension,
as well as the endpoint Sobolev embedding
$W^{1,1}\hookrightarrow L^\infty$.
Estimates relating the free energy to its dissipation were also used in \cite{Wi13,Wi20} to establish finite-time blowup in higher dimensions. In those works, the energy is estimated from below in terms of the dissipation, whereas here we derive an upper bound for the localized energy to establish uniform-in-time boundedness.
\end{remark}
To prove Proposition \ref{prop1}, we first establish
the following two lemmas.
\begin{lemma}\label{lem;interior-lyapunov}
	Let $(u,v)$ be the solution to \eqref{eqn;KS}.
	Then,
	\begin{align*}
	\frac{\mathrm{d}}{\mathrm{d} t}\mathcal F_R(t)+\mathcal{D}_R(t)
	=\,&-\int_{ \mathbb R^2}u(t,x)(\ln u(t,x) -v(t,x))\nabla(\ln u (t,x)- v(t,x))\cdot \nabla \varphi_R^2\,\mathrm{d}x
	\\
	&-\int_{ \mathbb R^2}\partial_t v (t,x)\nabla v(t,x) \cdot \nabla\varphi_R^2(x)\,\mathrm{d}x.
	\end{align*}
\end{lemma}
The proof of Lemma \ref{lem;interior-lyapunov} can be found in  \cite[Lemma 4.12]{Ho26} for instance.

\begin{lemma}\label{lem;log}
Let $(u,v)$ be the radially symmetric solution to \eqref{eqn;KS} with $\|u_0\|_1<8\pi$. Let $R>0$ satisfying
\begin{align}
A=A(R):= \sup_{t \ge t_0} u(t,R)<\infty.
\label{eqn;A}
\end{align}
Then, there exists a constant $C(A)=C(A,R,t_0)>0$ such that
\begin{align*}
\int_{B_R(0)} u(t,x)\ln u(t,x) \,\mathrm{d}x\le\, C(A)(1+ \sqrt{\mathcal{D}_R(t)})
\end{align*}
and
\begin{align*}
\int_{B_R(0)} |\nabla v(t,x)|^2\,\mathrm{d}x \le\, C(1+\mathcal{D}_R(t)).
\end{align*}
\end{lemma}
\begin{remark}
By Lemma \ref{lem;log}, there holds
\begin{align*}
\mathcal{F}_R(t) \le\, C\left(  1+\mathcal{D}_R(t)\right).
\end{align*}
\end{remark}
\begin{remark}
It is worthwhile to mention that $R$ can be chosen arbitrarily if \eqref{eqn;A} is satisfied.
Once we find a radius $R>0$ such that $A<\infty$, the inner regularity estimates on $B_R(0)$ will be controlled. 
\end{remark}
\begin{pr}{Lemma \ref{lem;log}}
We first define by $m(t,r)$, $n(t,r)$ and $z(t,r)$ the mass functions as
\begin{align*}
m(t,r):=2\pi\int_0^r u(t,s)s \,\mathrm{d}s,\qquad
n(t,r):=2\pi\int_0^r v(t,s)s\,\mathrm{d}s,\qquad
z(t,r):=2\pi \int_0^r \partial_tv(t,s)s \,\mathrm{d}s.
\end{align*}
Note that $m'(t,r)=2\pi u(t,r)r \ge 0$ and $0\le m(t,r)\le \|u_0\|_1$ as well as it follows from the second equation in \eqref{eqn;KS} that
\begin{align*}
\tau z(t,r)=2\pi v_r(t,r)r-\lambda n(t,r) +m(t,r),
\end{align*}
so that
\begin{align*}
v_r(t,r)=- \frac{m(t,r)-\tau z(t,r)-\lambda n(t,r)}{2\pi r}.
\end{align*}
Moreover,
\begin{align*}
|z(t,r)|=\left|\int_{B_r(0)} \partial_t v \,\mathrm{d}x\right|\le\, |B_r(0)|^{\frac12} \left(\int_{|x|<r}|\partial_tv|^2\,\mathrm{d}x\right)^{\frac12}\le\, \sqrt{\pi} r \| \partial_tv\|_{L^2(B_r(0))}.
\end{align*}
Then we obtain from integrations by parts
\begin{align*}
\mathcal{I}(t,r):=\,&\int_{B_r(0)} u(t,x) x\cdot \nabla(\ln u(t,x) -v(t,x))\,\mathrm{d}x
\\
=\,&\int_{\partial B_r(0)} u(t,x) (x\cdot \nu)\,\mathrm{d}S-2\int_{B_r(0)}u(t,x)\,\mathrm{d}x-2\pi\int_0^ru(t,s)v_s(t,s)s^2\,\mathrm{d}s
\\
=\,&2\pi u(t,r)r^2-2m(t,r)+\frac{1}{2\pi}\int_0^r m'(t,s) \left[ m(t,s)-\tau z(t,s) -\lambda n(t,s)\right]\,\mathrm{d}s
\\
=\,&r m'(t,r)-2m(t,r)+\frac{1}{4\pi}m^2(t,r)
\\
&-\frac{\tau}{2\pi}\int_0^rm'(t,s) z(t,s)\,\mathrm{d}s
-\frac{\lambda}{2\pi}\int_0^r m'(t,s)n(t,s)\,\mathrm{d}s.
\eqntag
\label{eqn;I}
\end{align*}
On the one hand,
\begin{align*}
\left|\mathcal{I}(t,r)\right| \le\,&\left(\int_{B_r(0)}  |x|^2 u(t,x)\,\mathrm{d}x\right)^{\frac12} \left( \int_{B_r(0)}u(t,x) \left|\nabla (\ln u(t,x) -v(t,x))\right|^2\,\mathrm{d}x\right)^{\frac12}
\\
\le\,&r \sqrt{m(t,r)X_1} \,\le r \sqrt{\|u_0\|_1 X_1},
\end{align*}
where 
\begin{align*}
X_1:= \int_{B_R(0)}u(t,x) \left|\nabla (\ln u(t,x) -v(t,x))\right|^2\,\mathrm{d}x.
\end{align*}
On the other hand,
\begin{align*}
\left|\int_0^r m'(t,s)n(t,s)\,\mathrm{d}s\right|\le\,&\sqrt{\pi}r \|u_0\|_1\|v\|_2,
\end{align*}
and
\begin{align*}
\left| \int_0^rm'(t,s) z(t,s)\,\mathrm{d}s \right|\le\,&\sqrt{\pi} r \|u_0\|_1 \|\partial_t v\|_{L^2(B_r(0))}\le\,\sqrt{\pi}r \|u_0\|_1 X_2,
\end{align*}
where
\begin{align*}
X_2:= \|\partial_t v\|_{L^2(B_R(0))},
\end{align*}
from which we deduce from \eqref{eqn;I} that
\begin{align*}
\left(2-\frac{\|u_0\|_1}{4\pi}\right)m(t,r)
\le\,&\left(2-\frac{m(t,r)}{4\pi}\right)m(t,r)
\\
=\,&r m'(t,r)-\frac{\tau}{2\pi}\int_0^r m'(t,s)z(t,s)\,\mathrm{d}s-\frac{\lambda}{2\pi}\int_0^r m'(t,s)n(t,s)\,\mathrm{d}s-\mathcal{I}(t,r)
\\
\le\,&r m'(t,r)+\frac{\tau r \|u_0\|_1}{2\sqrt{\pi}}X_2+\frac{\lambda r \|u_0\|_1 \|v\|_2}{2\sqrt{\pi}} +r \sqrt{\|u_0\|_1 X_1},
\end{align*}
where under the mass constraint $\|u_0\|_1<8\pi$
\begin{align*}
\delta:=\left(2-\frac{\|u_0\|_1}{4\pi}\right)>0.
\end{align*}
Integrating this over $(0,R)$ implies that
\begin{align*}
\delta \int_0^R \frac{m(t,r) }{r}\,\mathrm{d}r\le\,&\|u_0\|_1+\frac{\tau R \|u_0\|_1}{2\sqrt{\pi}}X_2+\frac{\lambda R \|u_0\|_1B_1}{2\sqrt{\pi}} +R \sqrt{\|u_0\|_1 X_1},
\end{align*}
where $B_1:= \sup_{t \ge 0}\|v(t)\|_2$ for $\lambda>0$.
By Tonelli's theorem,
\begin{align*}
\int_0^R \frac{m(t,r)}{r}\,\mathrm{d}r=\,&2\pi\int_0^R\frac1r \int_0^r u(t,s)s\,\mathrm{d}s\mathrm{d}r
\\
=\,&2\pi\int_0^R u(t,s)s\int_s^R \frac{1}{r}\,\mathrm{d}r\mathrm{d}s
\\
=\,&2\pi\int_0^R u(t,s)s\ln\left(\frac Rs\right)\,\mathrm{d}s
=\, \int_{B_R(0)}u(t,x)\ln\left(\frac{R}{|x|}\right)\,\mathrm{d}x,
\end{align*}
so that
\begin{align*}
\delta \int_{B_R(0)} u(t,x)\ln\left(\frac R{|x|}\right)\,\mathrm{d}x\le\,&\|u_0\|_1+\frac{\tau R \|u_0\|_1}{2\sqrt{\pi}}X_2+\frac{\lambda R \|u_0\|_1B_1}{2\sqrt{\pi}}+R \sqrt{\|u_0\|_1 X_1}.
\end{align*}
Next, we set
\begin{align*}
h(t,r):= \sqrt{u(t,r)},\quad g(t,r):=\frac{u_r(t,r)-u(t,r)v_r(t,r)}{\sqrt{u(t,r)}}.
\end{align*}
Then, we see 
\begin{align*}
h_r(t,r) + b(t,r) h(t,r) = \frac12 g(t,r),\qquad b(t,r):=\frac{m(t,r)-\tau z(t,r)-\lambda n(t,r)}{4\pi r},
\end{align*}
which solves that
\begin{align*}
h(t,r)=h(t,R) \exp\left(\int_r^R b(t,s)\,\mathrm{d}s\right) -\frac12\int_r^R g(t,s) \exp\left(\int_r^s b(t,\rho)\,\mathrm{d}\rho\right)\,\mathrm{d}s.
\end{align*}
Since the similar computations lead us to by use of negativity of $-\lambda n \le 0$ 
\begin{align*}
\int_r^R b(t,s)\,\mathrm{d}s\le\,&\frac{\|u_0\|_1}{4\pi}\ln\left(\frac R r\right)+\frac{\tau (R-r)}{4\sqrt{\pi}} X_2,
\end{align*}
and
\begin{align*}
\int_r^R |g(t,s)|\,\mathrm{d}s\le\,&\left(\int_r^R|g(t,s)|^2s\,\mathrm{d}s\right)^{\frac12} \left(\int_r^R \frac{1}{s}\,\mathrm{d}s\right)^{\frac12}\le\,\sqrt{\frac{X_1}{2\pi}} \sqrt{\ln\left(\frac Rr\right)},
\end{align*}
we have
\begin{align*}
\sqrt{u(t.r)}:=h(t,r)\le\,&\left[  \sqrt{u(t,R)} +  \sqrt{\frac{X_1}{8\pi}} \sqrt{\ln\left(\frac Rr\right)}\right] \left(\frac Rr\right)^{\frac{\|u_0\|_1}{4\pi}} \exp\left(\frac{\tau R }{4\sqrt{\pi}} X_2\right),
\end{align*}
which implies that
\begin{align*}
u(t,r)\le\,&\left[  \sqrt{u(t,R)} +  \sqrt{\frac{X_1}{8\pi}} \sqrt{\ln\left(\frac Rr\right)}\right]^2 \left(\frac Rr\right)^{\frac{\|u_0\|_1}{2\pi}} \exp\left(\frac{\tau R }{2\sqrt{\pi}} X_2\right)
\\
\le\,&C(A) (1+X_1)\left(1+\ln \left(\frac Rr\right)\right) \left(\frac Rr\right)^{\frac{\|u_0\|_1}{2\pi}}\exp\left(\frac{\tau R }{2\sqrt{\pi}} X_2\right),
\end{align*}
where $A$ is defined in \eqref{eqn;A}.
Hence,
\begin{align*}
\ln u(t,r) \le\,& \ln C(A) +\frac{\tau R}{2\sqrt{\pi}}X_2+\ln(1+X_1)+\ln\left(1+\ln\left(\frac R r\right)\right) +\frac{\|u_0\|_1}{2\pi}\ln \left(\frac{R}{r} \right),
\end{align*}
and accordingly
\begin{align*}
\int_{B_R(0)}u(t,x)\ln u(t,x)\,\mathrm{d}x\le\,&\|u_0\|_1 \left[   \ln C(A) +\frac{\tau R}{2\sqrt{\pi}}X_2+\ln(1+X_1) \right]
\\
&+\left(1+\frac{\|u_0\|_1}{2\pi}\right)\int_{B_R(0)}u(t,x) \ln\left(\frac R{|x|}\right)\,\mathrm{d}x
\\
\le\,&C(A,R) \left( 1+ \sqrt{X_1}+X_2 \right)\le\, C(A,R)\left(1+\sqrt{\mathcal{D}_R(t)}\right),
\end{align*}
where we infer from the definition of $\varphi_R$ 
\begin{align*}
X_1=\int_{B_R(0)}u(t,x) \left|\nabla (\ln u(t,x) -v(t,x))\right|^2\,\mathrm{d}x\le\, \int_{\mathbb{R}^2}u(t,x) \left|\nabla (\ln u(t,x) -v(t,x))\right|^2\varphi_R^2(x)\,\mathrm{d}x
\end{align*}
and 
\begin{align*}
X_2=\|\partial_t v\|_{L^2(B_R(0))} \le\, \|\partial_t v \varphi_R\|_2.
\end{align*}
Finally,
\begin{align*}
\int_{B_R(0)} |\nabla v(t,x)|^2\,\mathrm{d}x=\,&2\pi\int_0^R\left| \frac{m(t,r)-\tau z(t,r)-\lambda n(t,r)}{2\pi r} \right|^2 r \,\mathrm{d}r
\\
\le\,&\frac{1}{\pi}\int_0^R \frac{m^2(t,r)}{r}\,\mathrm{d}r + \frac{\tau^2}{\pi}\int_0^R \frac{z^2(t,r)}{r}\,\mathrm{d}r+\frac{\lambda^2}{\pi}\int_0^R \frac{n^2(t,r)}{r}\,\mathrm{d}r
\\
\le\,&\frac{\|u_0\|_1}{\pi}\int_{B_R(0)} u(t,x) \ln\left(\frac R{|x|}\right)\,\mathrm{d}x+\frac{\tau^2 R^2}{2} X_2^2+\frac{\lambda^2R^2}{2}B_1^2
\\
\le\,&C(R)\left( 1+ \mathcal{D}_R(t) \right),
\end{align*}
which ends the proof.
\end{pr}

\vspace{3mm}
We are ready to show the proof of Proposition \ref{prop1}.

\begin{pr}{Proposition \ref{prop1}}
By the assumptions on the boundedness of solutions, that is,
\begin{align*}
\sup_{t \ge t_0}\| \sqrt{u(t) } (\ln u(t)-v)\|_{L^2(R< |x|<2R )}+ \sup_{t \ge t_0}u(t,R)\le\, C(R)\qquad \text{for $t\ge t_0$,}
\end{align*}
where $C$ is independent of $t$, Lemma \ref{lem;interior-lyapunov} and Lemma \ref{lem;radial;heat-anul} imply that
\begin{align*}
\frac{\mathrm{d}}{\mathrm{d} t}\mathcal{F}_R(t)+\frac12\mathcal{D}_R(t)\le\, C(R,t_0),\qquad t\ge t_0>0
\end{align*}
since $\nabla\varphi_R^2=-2 \varphi_R\nabla\psi_{16R}$ from the definition of $\varphi_R$.
Moreover,  Lemma \ref{lem;log} along with the outer estimates give
\begin{align*}
\int_{ \mathbb{R}^2}u(t,x)\ln u(t,x) \varphi_R^2(x)\,\mathrm{d}x
=\,&\int_{ \mathbb R^2}u(t,x)\ln u(t,x) \chi_{B_R(0)}(x)\,\mathrm{d}x
\\
&+\int_{ \mathbb{R}^2}u(t,x)\ln u(t,x) \left(\varphi_R^2(x)-\chi_{B_R(0)}(x)\right)\,\mathrm{d}x
\\
=\,&\int_{ B_R(0)}u(t,x)\ln u(t,x)\,\mathrm{d}x
\\
&+\int_{R \le |x| \le 2R }u(t,x)\ln u(t,x) \left(\varphi_R^2(x)-\chi_{B_R(0)}(x)\right)\,\mathrm{d}x
\\
\le\,&C(R,A)(1+\mathcal{D}_R(t)),
\end{align*}
where $A$ is defined in \eqref{eqn;A}.
Similarly,
\begin{align*}
\int_{ \mathbb R^2}|\nabla v(t,x)|^2\varphi_R^2(x)\,\mathrm{d}x\le\,&C(R)(1+\mathcal{D}_R(t)).
\end{align*}
As a consequence, for any $\alpha>0$  
\begin{align*}
\frac{\mathrm{d}}{\mathrm{d} t}\mathcal{F}_R(t)+\alpha\mathcal{F}_R(t)+\frac12\mathcal{D}_R(t)\le\, C(R,t_0)+\alpha C(1+\mathcal{D}_R(t)),
\end{align*}
from which we infer by taking $\alpha>0$ sufficiently small such that $\beta:=1/2-\alpha C>0$,
\begin{align*}
\frac{\mathrm{d}}{\mathrm{d} t}\mathcal{F}_R(t)+\alpha\mathcal{F}_R(t)+\beta \mathcal{D}_R(t)\le\, C(R,t_0,\alpha) \qquad \text{for~~$t\ge t_0$.}
\end{align*}
Hence, we complete the proof.
\end{pr}

\subsection{Interior estimates for the localized entropy}
In this subsection, we combine the differential energy inequality
established in Proposition \ref{prop1} with a lower bound
for the localized free energy.
Under the subcritical mass condition \eqref{eqn;mass-constraint},
Trudinger--Moser's inequality allows us to control
the interior entropy and the Dirichlet energy in terms
of $\mathcal{F}_R$.
Together with the uniform-in-time upper bound for $\mathcal{F}_R$,
this yields uniform-in-time interior estimates and a bound
for the local energy dissipation integrated over time intervals
of unit length.
\begin{proposition}\label{prop2}
Suppose  that Proposition \ref{prop1} holds.
Assume  that $u_0$ satisfies \eqref{eqn;mass-constraint}.
Then, there exists $C=C(R,t_0)>0$ independent of $t$ such that
\begin{align*}
\sup_{t \ge t_0}| \mathcal{F}_R(t) | \le \, C.
\end{align*}
Consequently,
\begin{align*}
\int_{ \mathbb{R}^2} u(t,x) \ln u(t,x) \varphi_R^2(x)\,\mathrm{d}x+ \int_{t}^{t+1} \int_{ \mathbb{R}^2} |\partial_t v (s, x)|^2\varphi_R^2(x) \,\mathrm{d}x\mathrm{d}s\le\, C,\qquad t\ge t_0.
\end{align*}
\end{proposition}

\begin{pr}{Proposition \ref{prop2}}
Thanks to the consequence of Proposition \ref{prop1}, together with positivity of~$\mathcal{D}_R(t)$, we obtain
the uniform-in-time upper bound that
\begin{align*}
\mathcal{F}_R(t) \le\,&\mathcal{F}_R(t_0)+C, \qquad t \ge t_0.
\end{align*}
Next, we shall prove the lower bound on $\mathcal{F}_R(t)$.
The lower bound for $\mathcal{F}_R$ follows from
a standard energy argument, see \cite[Proposition 4.14]{Ho26} and also \cite{Mi13,Na-Se-Yo}.
We omit the details and obtain
\begin{align*}
\mathcal{F}_R(t)
\ge\,&\gamma\int_{\mathbb{R}^2}u(t,x)\ln u(t,x)\varphi_R^2(x)\,\mathrm{d}x+\frac12\left(1-\frac{\|u_0\|_1}{8\pi(1-\gamma)}\right)\int_{\mathbb{R}^2} |\nabla v(t,x)|^2\varphi_R^2(x)\,\mathrm{d}x-C(R).
\eqntag
\label{eqn;lower-bound}
\end{align*}
for some $C(R)>0$ and all  $\gamma>0$ satisfying
\begin{align*}
0<\gamma < \frac{8\pi-\|u_0\|_1}{8\pi},
\end{align*}
which ensures the positivity of the second term on the right hand side of \eqref{eqn;lower-bound}.
Using $x\ln x \ge -1/e$ for $x\ge0$,
we deduce from \eqref{eqn;lower-bound} along with the upper bound on $\mathcal{F}_R$ that
\begin{align*}
\sup_{t \ge t_0}| \mathcal{F}_R(t)| \le C(R,t_0).
\end{align*}
This thus implies that
\begin{align*}
\sup_{t \ge t_0}\int_{\mathbb{R}^2}u(t,x)\ln u(t,x)\varphi_R^2(x)\,\mathrm{d}x+\sup_{t \ge t_0} \int_{\mathbb{R}^2} |\nabla v(t,x)|^2\varphi_R^2(x)\,\mathrm{d}x\le\, C(R,t_0).
\end{align*}
Finally, it follows from the boundedness on $\mathcal{F}_R(t)$ that
\begin{align*}
\int_t^{t+1} \int_{\mathbb{R}^2} |\partial_t v(s,x)|^2\varphi_R^2(x)\,\mathrm{d}x\mathrm{d}s \le\, C,\qquad t\ge t_0.
\end{align*}
The proof is complete.
\end{pr}

\vspace{5mm}
\section{Uniform-in-time estimates on the plane}\label{sect:4}
In this section, we prove Theorem \ref{thm}.
Combining the outer and inner estimates established in
Sections \ref{sect;2} and \ref{sect;3}, respectively, we obtain
the following uniform-in-time bound for radially symmetric
solutions to \eqref{eqn;KS}:
\begin{equation}\label{eqn;uniform-bound}
\sup_{t \ge t_0}
\left[
\int_{\mathbb{R}^2}
(1+u(t,x))\ln(1+u(t,x))\,\mathrm{d}x
+
\int_t^{t+1}\|\partial_t v(s)\|_2^2\,\mathrm{d}s
\right]
<\infty.
\end{equation}
Starting from \eqref{eqn;uniform-bound}, we complete the proof
by a standard iteration argument, which we briefly present
for completeness.

\begin{pr}{Theorem \ref{thm}}
	The uniform-in-time bound \eqref{eqn;uniform-bound} indeed
	follows from Lemma \ref{lem;out-log} and Proposition \ref{prop2} along with Lemma \ref{lem;infty-u} by picking $R>1$ satisfying \eqref{eqn;assumption-R}.
	Once the uniform bound for the modified entropy is established,
	an argument similar to that used in the proof of
	Lemma \ref{lem;out-log} yields
	(cf. \cite{Ho-La25,Na-Se-Ya})
	\[
	\sup_{t \ge t_0}
	\int_t^{t+1}\int_{\mathbb R^2}
	\frac{|\nabla u(s,x)|^2}{1+u(s,x)}
	\,\mathrm{d}x\,\mathrm{d}s
	<\infty.
	\]
	Together with mass conservation, this estimate and a suitable
	version of Gagliardo--Nirenberg's inequality
	(see \cite{Ga-Za} and \cite[Lemma 2.3]{Na-Ya20})
	imply
	\begin{equation}\label{eqn;ulogu}
	\sup_{t \ge t_0}
	\int_t^{t+1}\int_{\mathbb R^2}u^2(s,x)
	\,\mathrm{d}x\,\mathrm{d}s
	<\infty.
	\end{equation}
	Combining the differential inequality for $\|u(t)\|_2^2$;
	\begin{align*}
	\frac{\mathrm{d}}{\mathrm{d}t}\|u(t)\|_2^2+\|\nabla u\|_2^2\le C\|\partial_t v\|_2^2\|u\|_2^2+C
	\end{align*}
	with \eqref{eqn;uniform-bound} and \eqref{eqn;ulogu},
	we apply the uniform Gronwall inequality to obtain
	\begin{align*}
	\sup_{t \ge t_0}\|u(t)\|_2<\infty.
	\end{align*}
This further yields
\begin{align*}
u\in L^\infty\bigl((t_0,\infty);L^p(\mathbb R^2)\bigr)
\qquad\text{for every }p\in(2,\infty).
\end{align*}
Therefore,
	parabolic smoothing estimates for the second equation in \eqref{eqn;KS} give
	\[
	\nabla v\in
	L^\infty\bigl((t_0,\infty)\times\mathbb R^2\bigr).
	\]
	Finally, we exploit the Moser iteration \cite{Al} to conclude that
	\[
	u\in L^\infty\bigl((t_0,\infty)\times\mathbb R^2\bigr),
	\]
	which completes the proof.
\end{pr}

\section*{Acknowledgments}\label{sec.Acknowledgments}
The work of the  author is partially supported by JSPS KAKENHI Grant Number 25K17274 (Early-Career Scientists), 25KJ0279 (Research Fellowship for Young Scientists).

\appendix
\section{Radial lemmas}\label{sect;appendix}
We collect some useful lemmas for radial functions in this section, which are exploited in Section~\ref{sect;2} in order to obtain the outer estimates for solutions to \eqref{eqn;KS}. 
\begin{lemma}\label{lem;radial-lemma}
	Let $f$ be a radially symmetric  function in $L^p(\mathbb{R}^2)$ such that $|f|^\frac q2\in \dot{H}^1(\mathbb{R}^2)$ with $1\le p,q <\infty$.
	Then,
	\begin{align*}
	|f(x) |\le\, \left( \frac{p+q}{2\pi q |x|} \right)^{\frac2{p+q}} \|f\|_p^{\frac{p}{p+q}} \| \nabla |f|^\frac q2\|_2^{\frac2{p+q}},\qquad x \in \mathbb R^2\setminus\{0\}.
		\end{align*}
\end{lemma}
\begin{pr}{Lemma \ref{lem;radial-lemma}}
	It suffices to show the desired estimate for $f \in C^\infty_0(\mathbb{R}^2)$.
	Since it follows that
	\begin{align*}
	|f(x)|^{\frac {p+q}2} = - \frac{p+q}{q}\int_{|x|}^\infty| f(r)|^\frac p2 ( |f(r)|^\frac q2)' \,\mathrm{d}r,
	\end{align*}
	we obtain
	\begin{align*}
	|f(x)|^{\frac {p+q}2}  \le\,&\frac{p+q}{q} \left( \int_{|x|}^\infty \frac{|f(r)|^{p}}{r}\,\mathrm{d}r \right)^{\frac12} 
	\left( \int_{|x|}^\infty | (|f(r)|^\frac q2)'|^2 r \,\mathrm{d}r \right)^\frac12
	\\
	\le\,&\frac{p+q}{2\pi q |x|} \| f\|_p^{\frac p2}  \| \nabla |f|^\frac q2\|_2,
	\end{align*}
	which ends the proof.
\end{pr}

\begin{lemma}\label{lem;radial2}
Let $f$ be a radially symmetric function in $L^1(\mathbb{R}^2)\cap \dot{H}^1(\mathbb{R}^2)$ and let $\psi_R$ be the cut-off function defined in \eqref{eqn;cutoff}. 
Then, for $\alpha>0$
\begin{align*}
|f(x)\psi_R(x)| \le\, \frac{1}{2\pi|x|} \sqrt{1+\frac1\alpha} \|f\|_1^\frac12 \left( \int_{\mathbb{R}^2}\frac{|\nabla (f \psi_R)|^2}{1+|f|}\,\mathrm{d}y \right)^\frac12+\alpha,\qquad x\in\mathbb{R}^2\setminus\{0\}.
\end{align*}
\end{lemma}

\begin{pr}{Lemma \ref{lem;radial2}}
Let $f \in C^\infty_0(\mathbb{R}^2)$ for simplicity. Set $F:=f \psi_R$. Note that $F$ is also a radial symmetric function.
Then,  for any $\alpha>0$
\begin{align*}
|F(x)| \le\,&(|F(x)|-\alpha)_+ + \alpha
\\
\le\,&\int_{|x|}^\infty \chi_{ \{|F|>\alpha\} } (r)| F'(r)|\,\mathrm{d}r+\alpha
\\
\le\,& \frac{1}{|x|} \left( \int_{|x|}^\infty  \chi_{ \{|F|>\alpha\} } (r) (1+|f(r)|) r \,\mathrm{d}r \right)^{\frac12} \left( \int_{|x|}^\infty \frac{|F'(r)|^2}{1+|f(r)|}r \,\mathrm{d}r\right)^\frac12+\alpha
\\
\le\,&\frac{1}{2\pi|x|} \sqrt{1+\frac1\alpha} \|f\|_1^\frac12 \left( \int_{\mathbb{R}^2}\frac{|\nabla (f\psi_R)|^2}{1+|f|}\,\mathrm{d}y \right)^\frac12+\alpha.
\end{align*}
This is the desired inequality.
\end{pr}

\vspace{3mm}
The next lemma provides a uniform-in-time $\dot H^1$ estimate on annular regions for radially symmetric solutions to linear parabolic equations, with radially symmetric initial data in $\dot H^1$ and radially symmetric forcing terms in $L^1$.
Although such solutions may exhibit logarithmic singularity both near the origin and at spatial infinity and hence need not be  bounded in $\dot H^1(\mathbb R^2)$, even in $\dot{H}^1(B_R(0))$, uniformly in time, their gradients can be controlled uniformly on annular regions.
This estimate will be used in Section \ref{sect;2} to derive the outer estimates near the boundary of a suitably chosen ball.
\begin{lemma}\label{lem;radial;heat-anul}
	Let $f\in L^\infty(0,\infty;L^1(\mathbb R^2))$ and $w_0\in\dot H^1(\mathbb R^2)$ be radially symmetric, and let $w$ be the radially symmetric solution to
	\begin{align*}
		\tau\partial_t w-\Delta w+\lambda w=f
	\qquad\text{in }(0,\infty)\times\mathbb R^2,
	\end{align*}
	where $\tau>0$ and $\lambda\ge0$.
	Then there exists a constant $C=C(\tau,\lambda)>0$, independent of $t$, $R$, $w_0$, and $f$, such that
	\begin{align*}
		\sup_{t\ge0}\|\nabla w(t)\|_{L^2(A_R)}
	\le
	\|\nabla w_0\|_{2}
	+
	C\|f\|_{L^\infty(0,\infty;L^1(\mathbb R^2))}
	\end{align*}
	for every $R>0$, where
	\begin{align*}
	A_R:=\{x\in\mathbb R^2; R<|x|<2R\}.
	\end{align*}
\end{lemma}
\begin{remark}
The claim in Lemma \ref{lem;radial;heat-anul} holds true even by replacing the annular regions $A_R$  by $A_{kR}$ for any~$k>0$.
\end{remark}
\begin{pr}{Lemma \ref{lem;radial;heat-anul}}
For simplicity, let $M:=  \|f\|_{ L^\infty(0,\infty; L^1(\mathbb{R}^2)) } $.
Owing to the integral formulation for $w$ and semigroup estimates, we infer that
\begin{align*}
\|\nabla w(t)\|_{L^2(A_R)} \le\,& \| \nabla e^{\tau^{-1}t(\Delta-\lambda) } w_0\|_{L^2(A_R)}
+\int_0^{\tau^{-1}t} \| \nabla e^{ s(\Delta-\lambda) } f(t-\tau s) \|_{L^2(A_R)}\,\mathrm{d}s
\\
\le\,&\| \nabla w_0\|_{2}
+\int_0^{\infty} \chi_{ [0,\tau^{-1}t] }(s) \| \nabla e^{ s\Delta } f(t-\tau s) \|_{L^2(A_R)}\,\mathrm{d}s
\\
\le\,&\| \nabla w_0\|_2
+\int_0^{R^2}  \chi_{ [0,\tau^{-1}t] }(s) \| \nabla e^{ s\Delta } f(t-\tau s) \|_{L^2(A_R)}\,\mathrm{d}s
\\&+\int_{R^2}^\infty\chi_{ [0,\tau^{-1}t] }(s)  |A_R|^{\frac12} \| \nabla e^{ s\Delta } f(t-\tau s) \|_{L^\infty(A_R)}\,\mathrm{d}s
\\
\le\,&\|\nabla w_0\|_2+\int_0^{R^2} \chi_{ [0,\tau^{-1}t] }(s) \| \nabla e^{ s\Delta } f(t-\tau s) \|_{L^2(A_R)}\,\mathrm{d}s
+CRM\int_{R^2}^\infty s^{-\frac32} \,\mathrm{d}s
\\
\le\,&\|\nabla w_0\|_2+\int_0^{R^2} \chi_{ [0,\tau^{-1}t] }(s) \| \nabla e^{ s\Delta } f(t-\tau s) \|_{L^2(A_R)}\,\mathrm{d}s
+CM.
\eqntag
\label{eqn;w-anul}
\end{align*}
Hence, it suffices to check the integral on $[0,R^2]$.
The heat kernel $G_s$ and $f$ are radial symmetry, so is $e^{s\Delta}f=G_s*f$.
Let $\mathbf{e}_1=(1,0)$ and $\mathbf{e}_\theta=(\cos\theta, \sin\theta)$ for $\theta \in [0,2\pi]$.
Since it follows from the polar coordinate transformation that for $r>0$
\begin{align*}
(e^{s\Delta}f(t-\tau s) )(r)=\,& \int_0^\infty \int_0^{2\pi} G_s( r \mathbf{e}_1- \rho \mathbf{e}_\theta) f(t-\tau s, \rho)\rho \,\mathrm{d}\theta \mathrm{d}\rho
\\
=\,&\int_0^\infty f(t-\tau s, \rho) \rho \int_0^{2\pi} G_s( r \mathbf{e}_1- \rho \mathbf{e}_\theta) \,\mathrm{d}\theta \mathrm{d}\rho,
\end{align*}
the derivative of this is represented as
\begin{align*}
\partial_r (e^{s\Delta}f (t-\tau s)) (r) =\,&2\pi \int_0^\infty f(t-\tau s,\rho)\rho \mathcal{K}(s,r,\rho)\,\mathrm{d}\rho,
\end{align*}
where
\begin{align*}
\mathcal{K} (s,r,\rho):= \frac{1}{2\pi} \int_0^{2\pi } \partial_r G_s( r \mathbf{e}_1- \rho \mathbf{e}_\theta) \,\mathrm{d}\theta.
\end{align*}
Then, Minkowski's inequality implies that
\begin{align*}
\| \nabla e^{ s\Delta } f(t-\tau s)\|_{L^2(A_R)} =\,&\left(2\pi\int_{R}^{2R} \left| 2\pi \int_0^\infty f(t-\tau s,\rho)\rho \mathcal{K}(s,r,\rho)\,\mathrm{d}\rho \right|^2 r \,\mathrm{d}r \right)^{\frac12}
\\
\le\,&2\pi\int_0^\infty |f(t-\tau s,\rho)| \rho \left( 2\pi\int_{R}^{2R} \left| \mathcal{K}(s,r,\rho) \right|^2 r\,\mathrm{d}r\right)^{\frac12} \,\mathrm{d}\rho.
\eqntag
\label{eqn;anulas}
\end{align*}
By the explicit formula for the heat kernel,
\begin{align*}
\left| \nabla G_s(z)\right| \le\,& \frac{|z|}{2\pi s^2} \exp\left( -\frac{|z|^2}{4s} \right)\le\, \frac{C}{s \sqrt{s}}\exp\left(-\frac{|z|^2}{8s}\right), \qquad z\in \mathbb{R}^2,
\end{align*}
we observe that
\begin{align*}
\left|\partial_r G_s( r \mathbf{e}_1- \rho \mathbf{e}_\theta)  \right|\le\,&\frac{C}{s\sqrt{s}}\exp\left(-\frac{(r-\rho)^2}{8s}\right) \exp\left(-\frac{r\rho}{4s}(1-\cos\theta) \right),
\end{align*}
where $|r \mathbf{e}_1-\rho\mathbf{e}_\theta|^2 = r^2+\rho^2-2r\rho\cos\theta= (r-\rho)^2+2r\rho(1-\cos\theta)$.
Hence, by using the symmetry of the integral and $1-\cos\theta = 2 \sin^2(\theta/2)\ge 2\theta^2/\pi^2$, we have
\begin{align*}
\int_0^{2\pi} \left| \partial_r G_s( r \mathbf{e}_1- \rho \mathbf{e}_\theta)  \right|\,\mathrm{d}\theta
\le\,&\frac{C}{s\sqrt{s}}\exp\left(-\frac{(r-\rho)^2}{8s}\right)  \int_0^\pi\exp\left(-\frac{r\rho}{4s}(1-\cos\theta) \right)\,\mathrm{d}\theta
\\
\le\,&\frac{C}{s\sqrt{s}}\exp\left(-\frac{(r-\rho)^2}{8s}\right)  \int_0^\pi\exp\left(-\frac{r\rho}{2\pi^2s}\theta^2\right)\,\mathrm{d}\theta.
\eqntag
\label{eqn;K}
\end{align*}
Next, if $R/2<\rho < 3R$,
we further estimate \eqref{eqn;K} as
\begin{align*}
\int_0^\pi\exp\left(-\frac{r\rho}{2\pi^2s}\theta^2\right)\,\mathrm{d}\theta\le\,& C\sqrt{\frac{s}{r\rho }}\le\,C \sqrt{\frac{s}{r R}},
\end{align*} 
so that
\begin{align*}
\left|\mathcal{K}(s,r,\rho) \right| \le\,&\frac{C}{s\sqrt{r R}}\exp\left(-\frac{(r-\rho)^2}{8s}\right)  \qquad \text{for}\quad \frac{R}{2} < \rho <3R.
\end{align*}
This implies that
\begin{align*}
2\pi\int_{R}^{2R} \left| \mathcal{K}(s,r,\rho) \right|^2 r\,\mathrm{d}r
\le\,&\frac{C}{s^2R} \int_{0}^{\infty}\exp\left(-\frac{(r-\rho)^2}{4s}\right) \,\mathrm{d}r\le\, \frac{C}{s^{\frac32}R}.
\end{align*}
Otherwise, that is, for either $\rho\le R/2$ or $3R\le \rho$,
\begin{align*}
\left|\mathcal{K}(s,r,\rho) \right| \le\,&\frac{C}{s\sqrt{s}}\exp\left(-\frac{(r-\rho)^2}{8s}\right)  \int_0^\pi\exp\left(-\frac{r\rho}{2\pi^2s}\theta^2\right)\,\mathrm{d}\theta
\\
\le\,&\frac{\pi C}{s\sqrt{s}}\exp\left(-\frac{(r-\rho)^2}{8s}\right).
\end{align*}
Observing that $|r-\rho|\ge  R/2$ for $R<r <2R$ and for either $\rho\le R/2$ or $3R\le \rho$, we realize that
\begin{align*}
2\pi\int_{R}^{2R} \left| \mathcal{K}(s,r,\rho) \right|^2 r\,\mathrm{d}r
\le\,&\frac{C}{s^3}\int_R^{2R} \exp\left(-\frac{(r-\rho)^2}{4s}\right) r \,\mathrm{d}r
\\
\le\,&\frac{C}{s^3}\int_R^{2R} \exp\left(-\frac{R^2}{16s}\right) r \,\mathrm{d}r\le\, \frac{C R^2}{s^3}\exp\left(-\frac{R^2}{16s}\right).
\end{align*}
Consequently,
\begin{align*}
\left(  2\pi\int_{R}^{2R} \left| \mathcal{K}(s,r,\rho) \right|^2 r\,\mathrm{d}r  \right)^{\frac12}
\le\,
\left\{
\begin{aligned}
&\frac{C}{s^{\frac34}\sqrt{R}},\quad &&\text{$\frac R2<\rho<3R$},
\\
&\frac{C R}{s^{\frac32}}\exp\left(-\frac{R^2}{32s}\right),\quad &&\text{otherwise}.
\end{aligned}
\right.
\end{align*}
Therefore, \eqref{eqn;anulas} is estimated as
\begin{align*}
\| \nabla e^{ s\Delta } f(t-\tau s)\|_{L^2(A_R)} \le\,&2\pi\int_0^\infty |f(t-\tau s,\rho) |\rho \left( 2\pi\int_{R}^{2R} \left| \mathcal{K}(s,r,\rho) \right|^2 r\,\mathrm{d}r\right)^{\frac12} \,\mathrm{d}\rho
\\
=\,&2\pi\int_{0}^{\frac R2} |f(t-\tau s,\rho) |\rho \left( 2\pi\int_{R}^{2R} \left| \mathcal{K}(s,r,\rho) \right|^2 r\,\mathrm{d}r\right)^{\frac12} \,\mathrm{d}\rho
\\
&+2\pi\int_{\frac R2}^{3R} |f(t-\tau s,\rho) |\rho \left( 2\pi\int_{R}^{2R} \left| \mathcal{K}(s,r,\rho) \right|^2 r\,\mathrm{d}r\right)^{\frac12} \,\mathrm{d}\rho
\\
&+2\pi\int_{3R}^\infty| f(t-\tau s,\rho)| \rho \left( 2\pi\int_{R}^{2R} \left| \mathcal{K}(s,r,\rho) \right|^2 r\,\mathrm{d}r\right)^{\frac12} \,\mathrm{d}\rho
\\
\le\,&M \left(  \frac{C}{s^{\frac34}\sqrt{R}} + \frac{C R}{s^{\frac32}} \exp\left(-\frac{R^2}{32s}\right)\right),
\end{align*}
and we end up with
\begin{align*}
\int_0^{R^2} \chi_{ [0,\tau^{-1}t] }(s) \| \nabla e^{ s\Delta } f(t-\tau s) \|_{L^2(A_R)}\,\mathrm{d}s
\le\,&C M \int_0^{R^2} \left[  \frac{1}{s^{\frac34}\sqrt{R}} + \frac{R}{s^{\frac32}} \exp\left(-\frac{R^2}{32s}\right)  \right]\,\mathrm{d}s \le\, CM.
\end{align*}
The proof thus follows from \eqref{eqn;w-anul} together  with the above estimate on $[0,R^2]$.
\end{pr}

\vspace{5mm}
\typeout{} %
\bibliographystyle{siam}
\bibliography{HL_KS_BDD_ref}

\end{document}